\documentclass{article}

\usepackage[nonatbib, preprint]{neurips_2019}
\usepackage[utf8]{inputenc} 
\usepackage[T1]{fontenc}    
\usepackage{hyperref}       
\usepackage{url}            
\usepackage{booktabs}       
\usepackage{amsfonts}       
\usepackage{nicefrac}       
\usepackage{microtype}      
\usepackage{amsmath}

\usepackage{tikz}
\usetikzlibrary{arrows}
\usetikzlibrary{positioning}
\usepackage{float}

\usepackage{amssymb}
\usepackage{bm}
\usepackage{todonotes}
\usepackage{lineno}

\usepackage{amsthm}
\usepackage{blindtext}

\newcommand{\ad}[0]{auto-diff}
\newcommand{\bp}[0]{backprop}
\newcommand{\rhs}[0]{r.h.s.}
\newcommand{\lhs}[0]{l.h.s.}

\theoremstyle{definition}
\newtheorem{definition}{Definition}[section]

\makeatletter
\renewcommand*\env@matrix[1][*\c@MaxMatrixCols c]{%
  \hskip -\arraycolsep
  \let\@ifnextchar\new@ifnextchar
  \array{#1}}
\makeatother

\title{QR and LQ Decomposition Matrix Backpropagation Algorithms for Square, Wide, and Deep - Real or Complex - Matrices and Their Software Implementation}
\author{
  Denisa A.O. Roberts\thanks{Corresponding author}\\
  AI SpaceTime\\
  NYC, NY \\
  \texttt{d.roberts@aispacetime.org} \\
  \And
  Lucas R. Roberts \\
  Amazon \\
  NYC, NY \\
  \texttt{rlucas@amazon.com} \\
}

\begin{document}

\maketitle

\begin{abstract}

  This article presents matrix backpropagation algorithms for the QR decomposition of matrices $A_{m, n}$, that are either square $(m = n)$, wide $(m < n)$, or deep $(m > n)$, with rank $k=min(m, n)$. Furthermore, we derive novel matrix backpropagation results for the pivoted (full-rank) QR decomposition and for LQ decomposition for deep input matrices. Differentiable QR decomposition offers a numerically stable, computationally efficient method to solve least squares problems frequently encountered in machine learning and computer vision. Other use cases such as graph learning and network compression are listed in the article. Software implementation across popular deep learning frameworks (PyTorch \cite{paszke2019pytorch}, TensorFlow \cite{abadi2016tensorflow}, MXNet \cite{chen2015mxnet}) incorporate the methods for general use within the deep learning community. Furthermore, this article aids the practitioner in understanding the matrix backpropagation methodology as part of larger computational graphs
  .
\end{abstract}

\section{Background}
\label{sec:intro}

 The QR decomposition is the thread that connects most of the algorithms of numerical linear algebra, including methods for least squares, eigenvalue, and singular value problems, as well as iterative methods for all of these and also for systems of equations ~\cite{trefethen1997numerical}. Despite the critical nature of the decomposition, the QR factorization and its gradient have lagged behind in deep learning research.  We surmise one reason is the absence of complete autodiff software implementations in the most common deep learning frameworks. The QR decomposition and the gradient of this decomposition have many uses in machine learning including for canonicalization of tensor networks in quantum machine learning (including the complex case) \cite{liao2019differentiable}, fitting least squares and Bayesian statistics \cite{seeger2017auto}, \cite{murphy2012machine}, and for optimum experimental design \cite{walter2010algorithmic}. Both QR (in the Householder implementation) and SVD solutions are backward stable and the one with the least computational cost can be chosen for full-rank least squares problems. The QR decomposition is a faster alternative \cite{trefethen1997numerical}, as well as a faster alternative to solving the least squares solution for approximately square matrices \cite{trefethen1997numerical}. Recent work \cite{9026903} compares SVD alternatives for solving least squares problems and provides an end to end solution to ellipse points fitting and human pose estimation applications and QR can be one. Other differentiable matrix factorizations were used in structured layers as part of end-to-end learning of computer vision models \cite{ionescu2015training}, or second-order pooling in graph neural networks \cite{ma2019graph}. The potential for other applications includes instances where least squares can approximate maximum likelihood estimation under assumptions relaxation, for practical reasons. Replacement of maximum likelihood estimation has been employed in some domains, for example with median rank regression \cite{olteanu2008technical}. Other potential use cases of differentiable QR decomposition are: complex uses cases in signal processing \cite{yasotharan2010simple}, network reconstruction \cite{li2017reconstruction}, neural networks for graph signal processing and graph clustering \cite{el2020orthonet}, Gaussian processes for relational reinforcement learning \cite{ramon2004numeric}, model compression in multitask learning \cite{kanakis2020reparameterizing}, transfer learning \cite{pan2009survey}.

In this document we refer to the process of calculating the gradient of the input matrix $\bm{A}$ interchangeably as QR matrix \bp\ or \ad\ QR. In \cite{seeger2017auto} an argument for the LQ \ad\ was presented for wide and square matrices. In \cite{liao2019differentiable} a different formula for the QR \ad\ was given. In \cite{walter2010algorithmic} the authors derive an analytical gradient for the square and deep real matrices case only. We fully derive new analytical formulae for all input matrix orders, clearly state the necessary assumptions, provide proofs, a supporting appendix, software implementations and numerical checks via central differences. 

The contributions of this article are:
 \begin{itemize}
    \item A review of the QR and LQ decompositions with a view toward \ad. We include a full derivation of the QR decomposition matrix \bp\ for square, wide, and deep matrices $\bm{A}$ (full-rank and pivoted full-rank), including their software implementations in popular deep learning frameworks.
    \item Using the partitioning trick, we derive a novel result for the LQ decomposition \bp\ for deep matrices. A \href{https://github.com/D-Roberts/lq_backprop}{Github} repository contains TensorFlow code for LQ \bp\ and corresponding numerical checks.
    \item A systematic process for deriving matrix \bp\ for the decomposition of deep, square, and wide inputs.
 \end{itemize}
 
 The exposition shows how the partitioning trick is applicable more broadly to cases where matrices are not of full rank. Once the input matrix is partitioned, the sub-matrices have full rank, thereby facilitating our derivations. The remainder of the article is structured as follows: in Section \ref{sec:prelim} we give an introduction to the QR and LQ decompositions. We also introduce the matrix \bp\ process and state basic results in matrix algebra that we use throughout the article. Then in Section 3 we derive expressions for the gradient through matrix backpropagation for QR (Deep and Square) case in Section 3.1, QR Wide case in Section 3.2, and LQ Deep case in Section 3.3. In Section 4 we give the complex QR matrix backprop. We then conclude with a few final notes. An extensive Appendix is available which includes derivations of results omitted for brevity.
 
\section{Preliminaries}\label{sec:prelim}
\label{one}
\subsection{The QR and The LQ Decomposition}
\label{prel}
The algebra of a QR decomposition is included in standard matrix algebra texts \cite{searle2017matrix}. 
Implementations of the QR decomposition are usually products of Householder rotations and are numerically stable. The implementation of the QR decomposition used in packages such as Numpy \cite{oliphant2006guide} and MXNet \cite{chen2015mxnet} typically wrap calls to the LAPACK  \cite{anderson1999lapack} library, as a sequence of two routines \textit{geqrf}, and \textit{orgqr}. The first LAPACK routine used in the sequence, $geqrf$, determines the Householder reflections whose matrix product determines $\bm{Q}$ and then returns the matrix $\bm{R}$ directly. The second LAPACK routine, \textit{orgqr}, returns the $\bm{Q}$ matrix. For GPU implementations, CuSolver (in the NVIDIA CUDA library \cite{nvidia2011nvidia}), is employed in the MXNet \cite{chen2015mxnet} implementation. Note that TensorFlow uses the \href{http://eigen.tuxfamily.org}{Eigen} library for linear algebra implementations. The decomposition typically returns, for an input matrix $\bm{A}$ of order $(m, n)$, matrices $\bm{Q}_{m, k}$ and $\bm{R}_{k, n}$ with $k=min(m, n)$, with $k$ typically the rank of the input matrix. This is the reduced mode decomposition. The first $k$ columns of Q form an orthonormal basis in the vector space spanned by the leading $k$ columns of the matrix $\bm{A}$. Notice that for wide input matrices (more columns than rows), $\bm{Q}$ is $m\times m$ with $k=m$, so $\bm{Q}$ is a square, orthogonal, full-rank matrix and all $m$ columns of $\bm{Q}$ are returned. In the case of deep matrices (more rows than columns, $m > n$), for the default (reduced) call to the decomposition, only the first $n$ columns of $\bm{Q}$ are included. In this case $\bm{R}$ is a square matrix while $\bm{Q}$ is not. This article assumes that the reduced mode QR (or LQ) is performed on the forward pass. The full (or complete) mode decomposition is also available and return matrices are of order $\bm{Q}_{m, m}$ and $\bm{R}_{m, n}$ but the reduced mode is typically the default because of computational efficiency considerations when solving least square problems \cite{trefethen1997numerical}. A further assumption we make throughout this article is that the rank of the input matrix is $k=min(m, n)$. Most of the article relies on the QR decomposition without column pivoting since, currently, the forward pass in the major deep learning frameworks do not have the pivoted algorithm implemented. However, such algorithms exist and the reader can consult for example \cite{van1983matrix} for theoretical details and the Eigen or SciPy \cite{virtanen2020scipy} libraries for software implementations.

There is a relationship between the QR and the LQ decomposition. If $\bm{A}=\bm{Q}\bm{R}$ is the QR decomposition of $\bm{A}$ then $\bm{A}^T = \bm{R}^T\bm{Q}^T$ is the LQ decomposition of $\bm{A}^T$.
The relationship implies that an algebraic result for the deep case of the LQ corresponds to the wide case of the QR and vice-versa. For ease of exposition and because of the aforementioned correspondence with the LQ decomposition, we focus our exposition primarily on the QR decomposition.

\subsection{Auto-differentiating Linear Algebra}
\label{ad}

\begin{figure}[!ht]
\centering
\begin{tikzpicture}[->,>=stealth',shorten >=1pt,auto,node distance=2cm,
  thick,main node/.style=
 {rectangle,
 draw=black, top color=white, bottom color=yellow!50,very thick, inner sep=0.5em, minimum size=3em, text centered}
  ]
  \tikzset{vertex/.style = {shape=rectangle,draw,minimum size=1.5em}}
  \tikzset{edge/.style = {->,> = latex'}}
  \node[main node, label=above:{$\bm{\bar{A}}\ \ $}] (A) {$\bm{A}$};
\node[above=0.3cm of A] (H) {};
  \node[main node, label=above:{\ \ $\bm{\bar{Q}}\ \ $}] (Q) [right of=A] {$\bm{Q}$};
\node[above=0.2cm of Q] (D) {};
  \node[main node, label=above:{\ \ $\bm{\bar{R}}\ \ $}] (R) [right of=Q] {$\bm{R}$};
 \node[above=0.2cm of R] (G) {};
  \node[main node] (F) [right of=R] {$\mathcal{L}$};
  \path[every node/.style={font=\sffamily\small,
  		fill=white,inner sep=3pt}]
    
    (F) edge [above right, bend right=75, dashed] node[right=5mm] {} (D)
    (F) edge [above right, bend right=60, dashed] node[right=10mm] {} (G)
    (G) edge [above, bend right=60, dashed] node[right=10mm] {} (H)
    (D) edge [above, bend right=50, dashed] node[left=10mm] {} (H);
    \draw[edge] (A) to[bend right=40] (Q);
    \draw[edge] (A) to[bend right=40] (R);
    \draw[edge] (R) to[bend right=40] (F);
    \draw[edge] (Q) to[bend right=40] (F);
\end{tikzpicture}
\caption{Forward pass calculations indicated with solid arrows and back-prop
calculations depicted with dashed arrows.
The matrices inside the nodes are the forward pass matrix decomposition values
while the matrices with bars are the \ad\ matrices.} 
\label{fig:adqr1}
\end{figure}
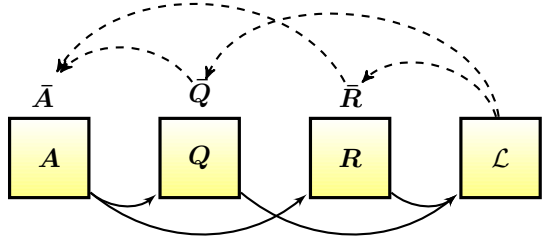

A collection of useful \ad\ results for linear algebra operators are given in \cite{seeger2017auto} and \cite{giles2008extended} and are used throughout this article. In \cite{ionescu2015training} a two step process for matrix \bp\ in the context of deep learning layers is described. We can view the calculation of the QR decomposition as successive computational layers. Each layer represents a matrix operation, with the matrix $\bm{A}$ preceding $\bm{Q}$ and $\bm{R}$ in the computational graph as depicted in Figure \ref{fig:adqr1}.


The goal of \ad\ through matrix \bp\ is to obtain an analytical formula for the gradient of matrix $\bm{A}$, given the gradient of matrices closer to the end of the topological order of the computational graph. In this article we refer to reverse mode \ad\ calculation of the gradient of $\bm{A}$, the QR matrix \bp\ process, as simply matrix \bp. From \cite{giles2008extended}, if $\bm{C}$ is an intermediate variable computed at a node in the computational graph and $\bm{C}$ is given as $\bm{C}=f(\bm{A}, \bm{B})$ a function of matrices $\bm{A}$ and $\bm{B}$ computed at a layer further ahead in the computational graph, then, from differential calculus,
\begin{equation} \label{eqn:dc}
    d\bm{C} = \frac{\partial f}{\partial \bm{A}}d\bm{A} + \frac{\partial f}{\partial \bm{B}}d\bm{B}. 
\end{equation}
In reverse mode \ad\ the infinitesimal perturbations are taken to be due to changes in the output $\mathcal{L}$. Note that $\mathcal{L}$ could be either a scalar valued loss function or an upstream layer in the computational graph. The sensitivities are computed starting at $\mathcal{L}$ and working backwards. By definition, an infinitesimal change in $\mathcal{L}$ and using the expression for $d\bm{C}$ in \ref{eqn:dc}, 
\begin{align*}
d\mathcal{L} 
&= Tr(\bar{\bm{C}}^Td\bm{C})\\
&= Tr(\bar{\bm{C}}^T\frac{\partial f}{\partial \bm{A}}d\bm{A}) + Tr(\bar{\bm{C}}^T\frac{\partial f}{\partial \bm{B}}d\bm{B}).
\end{align*}
We identify $\bar{\bm{A}} = \frac{\partial f}{\partial \bm{A}}^T\bar{\bm{C}}$ and $\bar{\bm{B}} = \frac{\partial f}{\partial \bm{B}}^T\bar{\bm{C}}$, which are the gradients sought.

Now assume that the gradient $\bar{\bm{C}}$ propagated from upstream in the topological ordering is the identity matrix $\bm{I}$. Then $\bar{\bm{A}} = \frac{\partial f}{\partial \bm{A}}^T$ and $\bar{\bm{B}} = \frac{\partial f}{\partial \bm{B}}^T$.
Since for the QR decomposition $\bm{A} = \bm{Q}\bm{R}$, with gradients of $\bm{Q}$ and $\bm{R}$ being backpropagated from upstream (denoted $\mathcal{L}$ in Figure \ref{fig:adqr1}), the tables are turned on the $\bm{C}=f(\bm{A},\bm{B})$ with $\bm{Q},\bm{R} = f(\bm{A})$.  The trace identity remains similar 
\begin{equation} \label{trid}
Tr(\bar{\bm{A}}^Td\bm{A}) =
    Tr(\bar{\bm{Q}}^Td\bm{Q}) + Tr(\bar{\bm{R}}^Td\bm{R}),  
\end{equation}
which is an identity we use repeatedly throughout the rest of the article.

The two-step \bp\ approach we repeatedly use in subsequent sections can be summarized as:
\begin{enumerate}
    \item Derive formulas for the variations of $\bm{Q}$ and $\bm{R}$ (or $\bm{L}$, $\bm{Q}$), denoted $d\bm{Q}$ and $d\bm{R}$, respectively, as a function of $d\bm{A}$. 
    \item Using the variations derived in step one, and the trace identity in Equation \ref{trid}, identify the gradient matrix, denoted $\bar{\bm{A}}$. 
\end{enumerate}

\subsection{A Collection of Useful Matrix Results}
\label{linalg}
Before proceeding to the QR and LQ matrix \bp\ derivations, it is useful to review a few matrix definitions and properties with proofs in any standard matrix algebra text such as \cite{searle2017matrix}. 

\begin{definition}\nonumber
If $\bm{A}$ is a square matrix, let
$sym(\bm{A}) = \frac{\bm{A} + \bm{A}^T}{2} $
denote a symmetric matrix.
\end{definition}
\begin{definition}\nonumber
 If a matrix $\bm{X}$ satisfies $\bm{X} = -\bm{X}^T$, then $\bm{X}$ is called a skew-symmetric matrix.
\end{definition}
Some useful properties of the trace operator:
\begin{enumerate}
\item $Tr(\bm{A}\bm{B}\bm{C}) = Tr(\bm{C}\bm{A}\bm{B}) = Tr(\bm{B}\bm{C}\bm{A})$, Invariance to Cyclic Permutations (ICP). 
\item $Tr(\bm{A}) = Tr(\bm{A}^T)$ Invariance to Transpose (IT).  
\item $Tr(\bm{A}(\bm{B}\circ \bm{C})) = Tr((\bm{B}^{\dagger} \circ \bm{A})\bm{C})$ general commutativity of Hadamard product (CH) for real and complex matrices.
\item From CH property and \cite{seeger2017auto}, for a square real matrix $\bm{M}$, $Tr(\bm{M}^T (sym(\bm{C}) \circ \bm{E})) = Tr(sym(\bm{M} \circ \bm{E})\bm{C})$, commutativity of products and symmetry (CPS). 
\end{enumerate}

Useful properties of matrices:
\begin{enumerate}
\item Similarly to \cite{seeger2017auto}, define a masking matrix $\bm{E}$ as
\begin{equation*}
e_{ij} = \begin{cases}0,  i < j\\
1,  i = j\\
2,  i > j 
\end{cases}
\end{equation*}
\item Similarly to \cite{seeger2017auto}, write $copyltu(\bm{A}) = sym(\bm{M} \circ \bm{E})$ (copy lower to upper) for a square real matrix $\bm{A}$. Note that $\bm{M} \circ \bm{E} = \bm{E} \circ \bm{M}$ as well.
\end{enumerate}
We will use these properties repeatedly in derivations in the sequel.

\section{QR and LQ Backpropagation Algorithms}

In this section we derive formulae for the gradient of $\bm{A}$ from the gradients of $\bm{Q}$ and $\bm{R}$ (or $\bm{L}$ and $\bm{Q}$ in the LQ decomposition) through matrix \bp . The matrix \bp\ algorithm takes as input the output of the forward pass of the decomposition ($\bm{Q}$, $\bm{R}$ or $\bm{L}$, $\bm{Q}$), in some cases the initial input matrix $\bm{A}$ (used in the partitioning trick), and the upstream gradients of $\bar{\bm{Q}}$, $\bar{\bm{R}}$ (or $\bar{\bm{L}}$, $\bar{\bm{Q}}$ in the LQ decomposition). 
We treat separately the QR decomposition \bp\ for deep (and square) matrices, the QR \bp\ for wide matrices, LQ \bp\ for wide (and square) matrices and finally LQ \bp\ for deep matrices. 
The derivation differs depending on the input matrix order (shape). As stated earlier, we assume rank $k$ of the input matrix.

Note that in all the derivations the forward implementation of the QR (or LQ) decomposition without pivoting is assumed. If the reduced pivoted decomposition is performed for input matrix $\bm{A}$ of full-rank $k=min\{m, n\}$, then $\bm{A}\bm{P} = \bm{Q}\bm{R}$, where $\bm{P}$ is a permutation matrix such that $\bm{Q}$ forms a basis for the first $k$ columns in $\bm{A}\bm{P}$, and $\bm{R}$ has the main diagonal elements in non-increasing order. In this case, another node is assumed to be included in the computational graph such that $\bm{A}\bm{P} = \bm{B}$ and then $\bm{B} = \bm{Q}\bm{R}$. To get the gradient of $\bm{A}$ on the backward pass, the derivations in the sequel are employed to get the gradient of matrix $\bm{B}$. Then, using the reverse mode \ad\ fundamentals in \cite{giles2008extended}, $\bar{\bm{A}} = \bar{\bm{B}} \bm{P}$. Note that the permutation matrix $\bm{P}$ does not receive a gradient from upstream computations.

The proceeding sections stand alone and can be independently consulted if the reader is interested in only one of the matrix orders or techniques presented. We dedicate the majority of the exposition to real matrices, and Section 4 demonstrates what changes for complex QR and LQ decomposition matrix backpropagation.

\subsection {QR Backpropagation: Real Square and Deep Matrices}\label{sec:qr_sq_deep}

In the forward pass, the QR decomposition of the input matrix A is obtained as described in Section \ref{sec:prelim}. On the backward pass the gradient is obtained using the two matrix \bp\ steps. A key assumption in the derivation for square and deep input matrix A is that the matrix  $\bm{R}$ is full rank. That $\bm{R}$ has full rank is a consequence of the assumption $\textnormal{rank}(\bm{A}) = \textnormal{min}(m,n)$, an assumption implicit as well in the derivations in \cite{seeger2017auto} for the wide input LQ decomposition gradient derivation.

\textbf{Proposition 1.}
Let $\bm{A}=\bm{Q}\bm{R}$ be the QR decomposition of a square or a deep matrix $\bm{A}$, $\bm{A} \in \mathbb{R}^{m, n}$, with $m \geq n$, such that $\bm{Q} \in \mathbb{R}^{m, n}$ and $\bm{R} \in \mathbb{R}^{n, n}$, with $\bm{R}$ an upper triangular matrix and $\bm{Q}$ orthogonal with $\bm{Q}^{T}\bm{Q} = \bm{I}_m.$ 

Then in the reverse mode \ad, 
\begin{equation}\label{eqn:a_bar}
\bar{\bm{A}} = \left[\bar{\bm{Q}} + \bm{Q}copyltu(\bm{M})\right]\bm{R}^{-T},    
\end{equation}
where $\bm{M} = \bm{R}\bar{\bm{R}}^T - \bar{\bm{Q}}^{T}\bm{Q}$. 

\textit{Proof.}
We apply the two step process described in Section ~\ref{ad}.
The proof is similar to the LQ decomposition \ad\ technique in \cite{seeger2017auto} for wide matrices. 
The LQ decomposition of a wide matrix parallels the QR decomposition of its transpose (a deep matrix).

\textbf{BP Step 1: Variations}
Calculate the variations $d\bm{Q}$ and $d\bm{R}$ for a given variation $d\bm{A}$.

\textit{Lemma 1.}
For a deep matrix $\bm{A}$ with QR decomposition $\bm{A}=\bm{Q}\bm{R}$, the variations of $\bm{Q}$ and $\bm{R}$ for a given variation $d\bm{A}$ of $\bm{A}$, $d\bm{Q}$ and $d\bm{R}$ respectively, are 
\begin{equation*}
\begin{split}
    d\bm{Q} &= (d\bm{A} - \bm{Q}d\bm{R})\bm{R}^{-1} \\
    d\bm{R} &= (sym(\bm{C})\circ \bm{E}^T)\bm{R}
\end{split}
\end{equation*}
with $\bm{C}=\bm{Q}^{T}d\bm{A}\bm{R}^{-1}.$

The proof of \textit{Lemma 1} is left to the Appendix. 

\textbf{BP Step 2: Partial Derivatives}
In the second step of the process for matrix \bp\ laid out in Section ~\ref{ad} one obtains partial derivatives using the variations from step one and the trace identity

\begin{equation}\label{eqn:tr_id}
Tr(\bar{\bm{A}}^Td\bm{A}) = Tr(\bar{\bm{Q}}^Td\bm{Q}) + Tr(\bar{\bm{R}}^Td\bm{R})
\end{equation}
The goal is to express $d\bm{Q}$ and $d\bm{R}$ variations on the \rhs\ of the trace identity as a function of $d\bm{A}$ and then identify the coefficients of $d\bm{A}$. First, replace the $d\bm{Q}$ expression in \textit{Lemma 1} on the \rhs\ to get
\begin{equation*}
 Tr(\bar{\bm{Q}}^Td\bm{Q}) + Tr(\bar{\bm{R}}^Td\bm{R}) = Tr(\bar{\bm{Q}}^T(d\bm{A} - \bm{Q}d\bm{R})\bm{R}^{-1}) + Tr(\bar{\bm{R}}^Td\bm{R}).
\end{equation*}
Next replace $d\bm{R}$ from \textit{Lemma 1}, use the ICP property of trace, and arrange the terms to isolate $d\bm{A}$. Considering the non $d\bm{A}$ terms only,
\begin{equation}
 Tr((\bar{\bm{R}}^T - \bm{R}^{-1}\bar{\bm{Q}}^T\bm{Q})(sym(\bm{C}) \circ \bm{E}^T)\bm{R})
    = %
    Tr(\bm{M}(sym(\bm{C}) \circ \bm{E}^T)),
\end{equation}

where we define $\bm{M} = \bm{R}\bar{\bm{R}}^T - \bar{\bm{Q}}^T\bm{Q}$. Next use the definition of $\bm{C}$ from \textit{Lemma 1}, as well as the ICP property from Section ~\ref{linalg} to simplify the trace to
\begin{equation*}
\begin{split}
Tr(\bm{M}(sym(\bm{C}) \circ \bm{E}^T)) & = Tr((sym(\bm{M}\circ \bm{E})^T\bm{C}) \\   
& = Tr((sym(\bm{M}\circ \bm{E})^T\bm{Q}^{T}d\bm{A}\bm{R}^{-1})\\
& = Tr(\bm{R}^{-1}(sym(\bm{M}\circ \bm{E})^T\bm{Q}^{T}d\bm{A}).
\end{split}
\end{equation*}
Returning to the \rhs\ expression in Equation \ref{eqn:tr_id}, containing both terms, we use ICP and IT to get
\begin{equation}\label{eqn:eq_rhs}
 Tr(\bm{R}^{-1}\bar{\bm{Q}}^Td\bm{A}) + Tr(\bm{R}^{-1}(sym(\bm{M}\circ \bm{E})^T\bm{Q}^{T}d\bm{A})= Tr([\bar{\bm{Q}} + \bm{Q}sym(\bm{M}\circ \bm{E})]\bm{R}^{-T})^Td\bm{A}).
\end{equation}
The matrix that left multiplies $d\bm{A}$ in  Equation ~\ref{eqn:tr_id} is $\bar{\bm{A}}^T$, and Equation \ref{eqn:eq_rhs} identifies
\begin{equation}
\begin{split}
\bar{\bm{A}} & = \left[\bar{\bm{Q}} + \bm{Q}sym(\bm{M}\circ \bm{E})\right]\bm{R}^{-T} \\
& = \left[\bar{\bm{Q}} + \bm{Q}copyltu(\bm{M})\right]\bm{R}^{-T}.
\end{split}
\end{equation}
Note that for real inputs the expression derived in Proposition 1 is equivalent to the result in \cite{walter2010algorithmic}. In \cite{walter2010algorithmic}, the QR gradient for deep or square real input matrices $\bm{A}$ is given in their Equation 42,
\begin{equation} \label{eqn:walter}
    \bar{\bm{A}} = \bm{Q}(\bar{\bm{R}} + \bm{P}_{L} \circ (\bm{R}\bar{\bm{R}}^T - \bar{\bm{R}}\bm{R}^T + \bm{Q}^T \bar{\bm{Q}} - \bar{\bm{Q}}^T \bm{Q})\bm{R}^{-T}) + (\bar{\bm{Q}} - \bm{Q} \bm{Q}^T \bar{\bm{Q}})\bm{R}^{-T}.
\end{equation}
Equation \ref{eqn:walter} (their Equation 42) simplifies to our Equation ~\ref{eqn:a_bar}. The matrix $\bm{P}_L = (i>j)$ in \ref{eqn:walter} is a strictly lower tridiagonal matrix with all ones below the diagonal and zeroes along and above the main diagonal. The equivalence of Equation \ref{eqn:a_bar} and Equation \ref{eqn:walter} for real matrices is proved in Appendix 6.2. The equivalence does not hold for complex inputs and Equation 42 needs a correction in the complex case.

\subsection {QR Backpropagation: Real Wide Matrices}

As in the previous section, we assume that the reduced mode QR decomposition is performed on the forward pass as described in Section ~\ref{prel}. On the backward pass, the gradient derivations for the wide case are complicated by the fact that $\bm{R}$ is no longer square and full rank. We employ a partitioning trick. Let $\bm{A}=\bm{Q}\bm{R}$ be the QR decomposition of a wide matrix $\bm{A}$, with $\bm{A} \in \mathbb{R}^{m, n}$, and $\bm{Q} \in \mathbb{R}^{m, m}$ a square orthogonal matrix, $\bm{R} \in \mathbb{R}^{m, n}$ an upper triangular matrix, and $m < n$. Throughout this article matrix $\bm{A}$ is assumed to be full rank $k=min(m, n)$. In the wide input case, a further assumption is that the first $k$ columns of A form a square full rank matrix. This is a strong assumption, generally met for random matrices. A pivoted QR decomposition, where $AP = QR$, can assure that the assumption is met. However the pivoted forward implementation is not available in the popular deep learning framework at the time of this writing. 

\textbf{Proposition 2.}

  On the backward pass let us assume that the input matrices $\bm{A}$, $\bm{R}$ and $\bm{\bar{R}}$ are partitioned as $\bm{A}=\left[\bm{X}|\bm{Y}\right]$, $\bm{R}=\left[\bm{U}|\bm{V}\right]$ and $\bm{\bar{R}}=\left[\bm{\bar{U}}|\bm{\bar{V}}\right]$ respectively with $\bm{X},\bm{R}, \bm{\bar{R}}\in \mathbb{R}^{m, m}$, and $\bm{Y}, \bm{V}, \bm{\bar{V}} \in \mathbb{R}^{m, n-m}$. Then $\bar{\bm{A}} = \left[\bar{\bm{X}} | \bar{\bm{Y}}\right]$, where
\begin{align} \label{eqn:prop2}
\bar{\bm{X}} &= [\bar{\bm{Q}}_{prime} + \bm{Q}copyltu(\bm{M})]\bm{U}^{-T}, \textnormal{   and}    \\
\bar{\bm{Y}} &= \bm{Q}\bar{\bm{V}},
\end{align}
with $\bar{\bm{Q}}_{prime} = \bar{\bm{Q}} + \bm{Y}\bar{\bm{V}}^T$, 
and $\bm{M} = \bm{U}\bar{\bm{U}}^T - \bar{\bm{Q}}_{prime}^{T}\bm{Q}$.

\textit{Proof}.

We follow the general process \bp\ outlined in Section ~\ref{ad} and employ the partitioning trick. For wide $\bm{A}$, $\bm{Q}$ is square, orthogonal, and full rank and the columns of $\bm{Q}$ form an orthogonal basis for the first $m$ columns of $\bm{A}$. We can formulate the QR decomposition as a multi-step sequential calculation which facilitates the gradient computation. The chain of operations is:
\begin{enumerate}
    \item Partition $\bm{A} = \left[\bm{X}|\bm{Y}\right]$, with $\bm{X}$ a square matrix of full rank $m$.
    \item Calculate the QR decomposition of $\bm{X}$, $\bm{X} = \bm{Q}\bm{U}$. The $\bm{Q}$ from the QR decomposition of $\bm{X}$ is the same $\bm{Q}$ in the QR decomposition of $\bm{A}$ produced on the forward pass. The matrix $\bm{U}$ is square and full-rank.
    \item Transform the remaining $n-m$ columns of $\bm{A}$, the $\bm{Y}$ partition, to get $\bm{V} = \bm{Q}^T\bm{Y}$.
    \item The reverse mode \ad\ employs this sequence in reverse order to get the partitioned $[\bar{\bm{X}}| \bar{\bm{Y}}]$ gradient. Then the full $\bar{\bm{A}}$ gradient is obtained by a simple concatenation node in the computational graph.
\end{enumerate}

\tikzstyle{line}  = [draw, -latex']

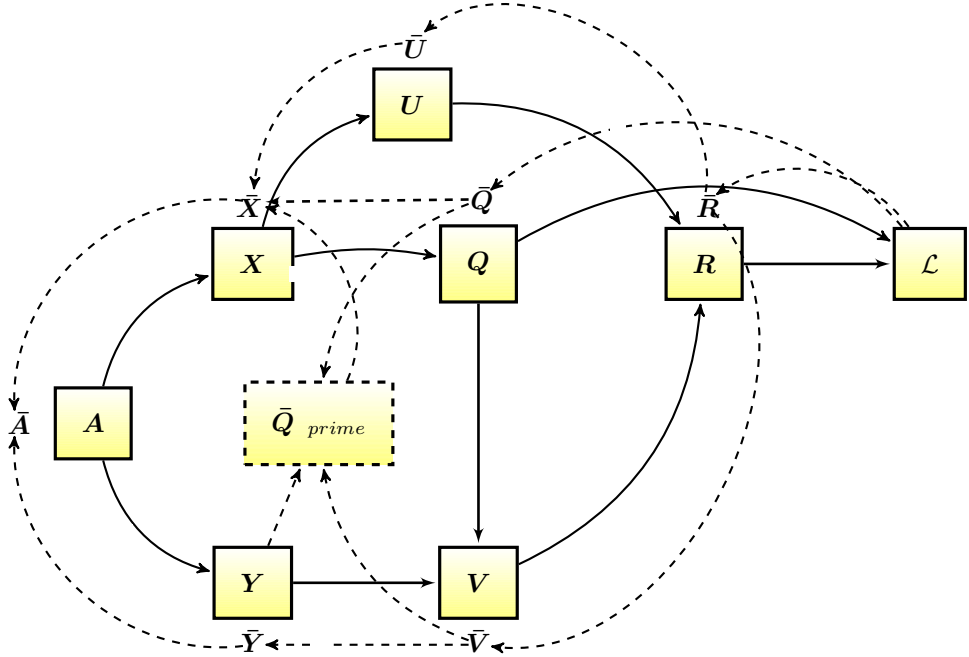
\begin{figure}
\centering
\begin{tikzpicture}[->,>=stealth',shorten >=1pt,auto,node distance=3cm,
  thick,main node/.style=
 {rectangle,
 draw=black, top color=white, bottom color=yellow!50, very thick, inner sep=1em, minimum size=1em, text centered}
  ]
  \node[main node, label=left:{$\bm{\bar{A}}\ \ $}] (A) {$\bm{A}$};
  \node[left=0.38cm of A] (H) {};
  \node[main node, label={$\bm{\bar{X}}\ $}] (X) [above right of=A] {$\bm{X}$};
  \node[main node, dashed] (Q_prime) [right of=A]{$\bm{\bar{Q}}_{\ \ prime}$};
  \node[above=0.2cm of X] (Xa) {};
  \node[main node, label=below:{\ \ $\bm{\bar{Y}}\ \ $}] (Y) [below right of=A] {$\bm{Y}$};
  \node[below=0.2cm of Y] (Yb) {};
  \node[main node, label=above:{\ \ $\bm{\bar{Q}} \ $}] (Q) [right of=X] {$\bm{Q}$};
  \node[above=0.2cm of Q] (Qa) {};
  \node[main node, label=above:{\ \  $\bm{\bar{U}}\  $}] (U) [above right of=X] {$\bm{U}$};
  \node[above=0.2cm of U] (Ua) {};  
  \node[main node, label=below:{\ \ $\bm{\bar{V}}\ \ $}] (V) [right of=Y] {$\bm{V}$};
  \node[below=0.2cm of V] (Vb) {};
  \node[main node, label=above:{\ \  $\bm{\bar{R}}\  $}] (R) [right of=Q] {$\bm{R}$};
  \node[above=0.2cm of R] (Ra) {};
  \node[main node] (F2) [right of=R] {$\mathcal{L}$};
  
  \path[every node/.style={font=\sffamily\small,
  		fill=white,inner sep=3pt}]

    (A) edge [bend right=30] node[left=1mm] {} (Y)
    (A) edge [bend left=30] node[left=10mm] {} (X)
    (X) edge [bend left=30] node[right=5mm] {} (U)
    (X) edge [above right, bend left=10, solid] node[] {} (Q)
    (V) edge [bend right=30] node[left=10mm] {} (R)
    (U) edge [bend left=30] node[left=5mm] {} (R)
    (Q) edge [bend left=30, solid] node[right=12mm] {} (F2)
        
    (Yb) edge [above, bend left=50, dashed] node[left=5mm] {} (H)
    (Xa) edge [above, bend right=50, dashed] node[left=5mm] {} (H)
    (Vb) edge [below, bend right=0, dashed] node[left=5mm] {} (Yb)
    (Ua) edge [left, bend right=40, dashed] node[left=5mm] {} (Xa)
    (F2) edge [left, bend right=45, dashed] node[left=10mm] {} (Ra)
    (F2) edge [right, bend right=45, dashed] node[left=10mm] {} (Qa)
    (Vb) edge [left, bend left=30, dashed] node[left=5mm] {} (Q_prime)
    (Qa) edge [left, bend right=30, dashed] node[left=5mm] {} (Q_prime)
    (Y) edge [above, dashed] node[left=5mm] {} (Q_prime)
    (Q_prime) edge [left, bend right=50, dashed] node[left=5mm] {} (Xa)
    (Ra) edge [right, bend left=70, dashed] node[above right=15mm] {} (Vb)
    (Ra) edge [left, bend right=70, dashed] node[right=5mm] {} (Ua)
    ;
     \path [line, line width=0.35mm] (Y)         -- (V);
     \path [line, line width=0.35mm] (R)         -- (F2);
     \path [line, line width=0.35mm] (Q)         -- (V);
     
     \path [line, line width=0.35mm, dashed ] (Qa)         -- (Xa);
\end{tikzpicture}
\caption{A graphical depiction of the QR auto-diff computation. 
Forward pass results are indicated with solid arrows and back-prop
calculations depicted with dashed arrows. Note that the partitioning is only employed as a device on the backward pass. The matrices inside the nodes are the forward pass matrix decomposition values
while the matrices with overbars are the gradient matrices. In the backward pass we employ the split and concatenate operations: $\bm{A} = [ \bm{X} | \bm{Y} ]$, $\bm{R} = [ \bm{U} | \bm{V} ]$, $\bm{\bar{R}} = [ \bm{\bar{U}} | \bm{\bar{V}} ]$ to finally obtain $\bm{\bar{A}} = [ \bm{\bar{X}} | \bm{\bar{Y}} ]$.}

\label{fig:wide_case_dag}
\end{figure}
The reader may find Figure \ref{fig:wide_case_dag} instructive to study for the \bp\ computations with a partitioning trick. Next we apply the chain rule and the two-step \bp\ algorithm to obtain $\bar{\bm{A}}$.

\textbf{BP Step 1: Variations} 
First calculate the variations of $\bm{Q}$ and $\bm{R}$, $d\bm{Q}$ and $d\bm{R}$, for a given variation $d\bm{A}$ of $\bm{A}$. 

\textit{Lemma 2.}
\label{lem:variations}
The variations of $\bm{A}$ and $\bm{R}$ are partitioned $d\bm{A}=[d\bm{X}|d\bm{Y}]$ and $d\bm{R} = [d\bm{U}|d\bm{V}]$. Then 
\begin{align}
    d\bm{Q} &= (d\bm{X} - \bm{Q}d\bm{U})\bm{U}^{-1} \\
    d\bm{V} &= \bm{Q}^{T}(d\bm{Y} - d\bm{Q}\bm{V}) \\
    d\bm{U} &= (sym(\bm{C})\circ \bm{E}^T)\bm{U},
\end{align}
with $\bm{C}=\bm{Q}^{T}d\bm{X}\bm{U}^{-1}$ and $\bm{Q}^Td\bm{Q}$ is skew-symmetric. 

\begin{proof}
The proof of \textit{Lemma 2} is left to the appendix. 
\end{proof}

\textbf{BP Step 2: Partial Derivatives}
On the second step of the matrix \bp\ process we derive $[\bar{\bm{X}}|\bar{\bm{Y}}]$ and consequently $\bar{\bm{A}}$. We begin by using the partitioned form to rewrite Equation \ref{eqn:tr_id},

\begin{equation}\label{eq:21}
Tr(\bar{\bm{X}}^Td\bm{X}) + Tr(\bar{\bm{Y}}^Td\bm{Y}) = Tr(\bar{\bm{Q}}^Td\bm{Q}) + Tr(\bar{\bm{U}}^Td\bm{U}) + Tr(\bar{\bm{V}}^Td\bm{V}).
\end{equation}

Replace the variations $d\bm{Q}$, $d\bm{U}$ and $d\bm{V}$ in the \rhs\ of Equation \ref{eq:21} with their expressions as functions of the $d\bm{X}$, and $d\bm{Y}$ from \textit{Lemma 2}. Now both the \lhs\ and \rhs\ are functions of only $d\bm{X}$ and $d\bm{Y}$ and are orthogonal to each other due to the nature of the partitioning. Matching the pre-multiplying matrices from the \lhs\ and \rhs\ of Equation \ref{eq:21} we identify the gradients $\bar{\bm{X}}$ and $\bar{\bm{Y}}$. 

Focusing on the \rhs of Equation \ref{eq:21}, first replace $d\bm{V}$ from \textit{Lemma 2} to get
\begin{equation*}
Tr(\bar{\bm{Q}}^Td\bm{Q}) + Tr(\bar{\bm{U}}^Td\bm{U}) - Tr(\bar{\bm{V}}^{T}\bm{Q}^{T}d\bm{Q}\bm{V}) + Tr(\bar{\bm{V}}^{T}\bm{Q}^{T}d\bm{Y}).
\end{equation*}

We identify the transpose of the $d\bm{Y}$ coefficient in Equation \ref{eq:21} as $\bar{\bm{Y}} = \bm{Q}\bar{\bm{V}}$. Now that we have identified the gradient $\bar{\bm{Y}}$, we omit the corresponding terms from both the \lhs\ and \rhs\ of Equation \ref{eq:21} and focus on identifying $d\bm{\bar{X}}$. We have
\begin{equation*}
Tr(\bar{\bm{X}}^Td\bm{X}) = Tr(\bar{\bm{Q}}^Td\bm{Q}) + Tr(\bar{\bm{U}}^Td\bm{U}) - Tr(\bar{\bm{V}}^{T}\bm{Q}^{T}d\bm{Q}\bm{V}).
\end{equation*}

The last trace can be simplified by replacing $\bm{Q}^Td\bm{Q}$, using skew-symmetry and $\bm{V} = \bm{Q}^{T}\bm{Y}$ from the partitioning argument, and the orthogonality of $\bm{Q}$ to get
\begin{align*}
Tr(\bar{\bm{V}}^{T}\bm{Q}^{T}d\bm{Q}\bm{V}) &= -Tr(\bar{\bm{V}}^{T}(d\bm{Q}^T\bm{Q})\bm{V}) \\
&= -Tr(\bar{\bm{V}}^{T}((d\bm{Q}^T\bm{Q})(\bm{Q}^{T}\bm{Y})) \\ 
&= - Tr(\bar{\bm{V}}^{T}d\bm{Q}^T\bm{Y}). 
\end{align*}

Further use the ICP and the IT property from Section \ref{linalg} to write 
\begin{equation*}
Tr(\bar{\bm{V}}^{T}d\bm{Q}^T\bm{Y}) = Tr(\bm{Y}\bar{\bm{V}}^{T}d\bm{Q}^T) = Tr(\bar{\bm{V}}\bm{Y}^Td\bm{Q}). 
\end{equation*}

With this result simplify the \rhs\ of Equation ~\ref{eq:21},
\begin{align} \label{eqn:xbar}
Tr(\bar{\bm{Q}}^Td\bm{Q}) + Tr(\bar{\bm{U}}^Td\bm{U}) + Tr(\bar{\bm{V}}\bm{Y}^Td\bm{Q})
&= Tr((\bar{\bm{Q}}^T + \bar{\bm{V}}\bm{Y}^T)d\bm{Q}) + Tr(\bar{\bm{U}}^Td\bm{U})\\
& = Tr(\bar{\bm{Q}}^T_{prime}d\bm{Q}) + Tr(\bar{\bm{U}}^Td\bm{U}),
\end{align}
where we use the notation $\bar{\bm{Q}}_{prime} \equiv \bar{\bm{Q}} + \bm{Y}\bar{\bm{V}}^T$. The rest of the argument for the $d\bm{X}$ term is similar to the argument in Section \ref{sec:qr_sq_deep} since $\bm{X}$ and $\bm{U}$ are square. Next, we follow the same process as in the proof of Proposition 1, with $\bm{X}$ in place of $\bm{A}$, $\bm{U}$ in place of $\bm{R}$ and $\bar{\bm{Q}}_{prime}$ in place of $\bar{\bm{Q}}$. We use $\bar{\bm{Q}}_{prime}$ instead of simply $\bar{\bm{Q}}$ since $\bm{Q}$ is used in multiple node computations in the computational graph showed in Figure \ref{fig:wide_case_dag}. Consequently, using the result in Section 3.1, we get $\bar{\bm{X}} = (\bar{\bm{Q}}_{prime} + \bm{Q} copyltu(\bm{M}))\bm{U}^{-T}$, with $\bm{M} = \bm{U}\bar{\bm{U}}^T - \bar{\bm{Q}}_{prime}^{T}\bm{Q}$. Now that we have identified both $\bar{\bm{X}}$ and $\bar{\bm{Y}}$, we concatenate them to obtain $\bar{\bm{A}}$
\begin{equation*}
\bar{\bm{A}} = \left[(\bar{\bm{Q}}_{prime} + \bm{Q} copyltu(\bm{M}))\bm{U}^{-T}|\bm{Q}\bar{\bm{V}}\right], 
\end{equation*}
as proposed. Notice that the calculation of $\bar{\bm{X}}$ for real matrices in the wide case employs the techniques in Section 3.1 hence either Equation \ref{eqn:a_bar} or the equivalent Equation \ref{eqn:walter} can be employed. In Appendix 5.2 we prove that Equation \ref{eqn:a_bar} and Equation \ref{eqn:walter} are equivalent for real matrices. However they are not equivalent for complex inputs and we will demonstrate in Section 4. Where possible we prefer to use Equation \ref{eqn:a_bar} for computational efficiency considerations.

The connection to the square case simplifies software implementations considerably. Although memory is re-used as much as possible for efficiency, additional memory must still be allocated for intermediary results. In the MXNet software implementation (C++) Equation \ref{eqn:a_bar} is implemented as a helper and then used directly to either obtain $\bar{\bm{A}}$ directly when $\bm{A}$ is square or deep, or to obtain $\bar{\bm{X}}$ when $\bm{A}$ is deep. TensorFlow had already implemented a QR backward method decomposition for square and deep input matrices (Python) and the Equation \ref{eqn:walter} was used. We implemented Equation \ref{eqn:walter} as a helper and calculated $\bar{\bm{X}}$ directly using the helper function. However, Equation 
\ref{eqn:a_bar} has a more computationally efficient implementation and should be preferred when possible. After obtaining $\bar{\bm{X}}$ and $\bar{\bm{Y}}$ we re-assemble the gradient $\bar{\bm{A}}$ via matrix concatenation. Numerical tests to test the correctness of the gradient via central differences were implemented as part of the typical test driven development approach. These tests can be referred to in the backend of the aforementioned deep learning frameworks on \href{https://github.com/tensorflow/tensorflow/blob/master/tensorflow/python/kernel_tests/qr_op_test.py}{Github}. 

\subsection{LQ Backpropagation: Real Square, Wide and Deep}

In \cite{seeger2017auto}, the reverse mode \ad\ for the LQ decomposition is derived for square and wide matrices, $\bm{A} \in \mathbb{R}^{m, n}$, $ m \le n$, as:
\begin{equation}
\bar{\bm{A}} = \bm{L}^{-T}\left[\bar{\bm{Q}} + copyltu(\bm{M})\bm{Q}\right],    
\end{equation}
where $\bm{M} = \bm{L}^T\bar{\bm{L}} - \bar{\bm{Q}}\bm{Q}^T$. A key assumption in the proof is that $\bm{L}$ is full rank. For the deep case, the reverse mode AD derivation is stymied by the shapes of $\bm{L}$ and $\bm{Q}$ and by the fact that $\bm{L}$ is not full rank. We again employ the partitioning trick to get the LQ decomposition gradient for deep input matrices. To the best of our knowledge this result is novel.

\textbf{Proposition 3}\label{sec:LQ}

Let $\bm{A}=\bm{L}\bm{Q}$ be the LQ decomposition of a deep matrix $\bm{A}$, with $\bm{A} \in \mathbb{R}^{m, n}$, $\bm{Q} \in \mathbb{R}^{n, n}$ a square orthogonal matrix, $\bm{L} \in \mathbb{R}^{m, n}$ a lower triangular matrix, and $m > n$, with the top $n$ rows of $\bm{A}$ forming a square full-rank matrix. In the reverse mode \ad\ $\bm{A}$, $\bm{L}$ and $\bm{\bar{L}}$ are partitioned into
$\bm{A} = \begin{bmatrix}
\bm{X} \\
\bm{Y}
\end{bmatrix}$, $\bm{L} = \begin{bmatrix}
\bm{U} \\
\bm{V}
\end{bmatrix}$, and $\bm{\bar{L}} = \begin{bmatrix}
\bm{\bar{U}} \\
\bm{\bar{V}}
\end{bmatrix}$, respectively
 with $\bm{X}, \bm{U}, \bm{\bar{U}}  \in \mathbb{R}^{n, n}$, $\bm{Y}, \bm{V}, \bm{\bar{V}} \in \mathbb{R}^{m-n, n}$.
 Then we calculate $\bar{\bm{A}}$ as the concatenation of $\bar{\bm{X}}$ and $\bar{\bm{Y}}$ along the row axis
 with
$\bar{\bm{X}}
= \bm{U}^{-T}[\bar{\bm{Q}}_{prime} + copyltu(\bm{M})\bm{Q}] $ and 
$\bar{\bm{Y}}
= \bar{\bm{V}}\bm{Q}$,  
with 
$\bm{M} = \bm{U}^T\bar{\bm{U}} - \bar{\bm{Q}}_{prime}\bm{Q}^T$, and  $\bar{\bm{Q}}_{prime} = \bar{\bm{Q}} + \bm{\bar{V}^T}\bm{Y}$. 

The proof of Proposition 3 parallels the proof of Proposition 2 and is left to the Appendix. The LQ decomposition and its gradient can be obtained using the TensorFlow based code in this Github \href{https://github.com/D-Roberts/lq_backprop}{repository}. 

\section{Complex Input QR Matrix Backpropagation}

In this section we derive the matrix backpropagation process for the reduced mode QR decomposition of full-rank (row or column) complex matrices of any shape (matrix order). 
Whether it is row or column rank will be inferred. 
The section is meant to serve as a stand-alone complete reference for the complex QR matrix backprop. 
The LQ derivations follow a similar reasoning pattern and will not be written explicitly. 
In this section we first briefly review the matrix backprop process and how the complex case differs from the real case.
Then we outline a few complex differentials fundamental concepts and
then derive the correct square and deep matrix QR reverse mode gradient when A is a complex input matrix.
We give a correction term to be added to the real deep case expression implemented in TensorFlow and PyTorch from \cite{walter2010algorithmic}. 
This correction is necessary for implementation to produce the correct complex gradient if the \cite{walter2010algorithmic} formula is used to calculate the gradient. 
We then outline the wide complex matrix backprop.

\subsection{General Matrix Backprop Process with Complex Inputs}
The general matrix backpropagation process remains as outlined in the real case, however there are a number of differences in the derivations of the reverse mode auto-diff for complex QR.
On the forward pass, the QR decomposition of complex input matrix returns the Q matrix with $\bm{Q}^{\dagger}\bm{Q}=\bm{I}$ and the upper triangular complex R matrix which has real main diagonal elements.

Recall from Section 2.2 that our goal in matrix backpropagation is to derive an analytical expression for the reverse mode gradient $\bar{\bm{A}}$. 
As in the real case, we assume a computational graph that optimizes upstream a real scalar loss function $\mathcal{L}$. Assuming the computational graph for the real case in Figure 1, $\mathcal{L}$ is a function of the input matrix $\bm{A}$ as well as the matrices $\bm{Q}$ and $\bm{R}$ at different layers of the graph. We can express an infinitesimal perturbation in the loss function $\mathcal{L}$ as 
\begin{equation*} 
    d\bm{\mathcal{L}} = \frac{\partial f}{\partial \bm{A}}d\bm{A} = Tr(\bar{\bm{A}}^{T}d\bm{A})
\end{equation*}

Since $\bm{A}=\bm{Q}\bm{R}$, $\bm{Q}$ and $\bm{R}$'s gradients are backpropagating from $d\mathcal{L}$, in the real case 
\begin{equation}\label{trid1}
Tr(\bar{\bm{A}}^{T}d\bm{A}) =
    Tr(\bar{\bm{Q}}^{T}d\bm{Q}) + Tr(\bar{\bm{R}}^{T}d\bm{R}),  
\end{equation}
then we identify the gradient of interest, $\bar{A}$, in the L.H.S term of Equation \ref{trid1} as the pre-multiplier of $d\bm{A}$. To do this identification we follow the standard matrix backprop two steps:

\textbf{Step 1: Variations}. Express the variations $d\bm{R}$ and $d\bm{Q}$ as functions of $d\bm{A}$.

\textbf{Step 2: Partial Derivatives}. Using the variations derived in step one, and the trace identity in Equation \ref{trid1}, identify the gradient matrix, denoted $\bar{\bm{A}}$. 
 
We now outline the differences from the real case. Firstly for the complex case the trace identity in Equation \ref{trid1} becomes:
\begin{equation}\label{trid2}
Tr(\bar{\bm{A}}^{\dagger}d\bm{A} + ({\bar{\bm{A}}}^{\dagger}d\bm{A})^{\dagger}) =
    Tr(\bar{\bm{Q}}^{\dagger}d\bm{Q} + h.c.) + Tr(\bar{\bm{R}}^{\dagger}d\bm{R} + h.c.).  
\end{equation}
For brevity, we use the notation $h.c.$ to mean that all the Hermitian conjugate terms of the existing terms in the parenthesis are included as well. For example $\bm{A} + \bm{B} + h.c. = \bm{A} + \bm{B} + \bm{A}^{\dagger} + \bm{B}^{\dagger}$. Our goal is the same as in the real case - to determine an expression for $\bar{\bm{A}}$. Note that $\bar{\bm{A}}$ is the pre-multiplier of $d\bm{A}$. One could work with the h.c. terms instead and get the expression for $\bar{\bm{A}}^{\dagger}$. But why do we have the extra $h.c.$ terms to begin with? The next section sheds some light on the inclusion of these additional terms. 

\subsection{Complex Differentiability Fundamentals}

We outline a few complex matrices and complex differentiability  concepts necessary in the derivation of complex QR matrix backprop. One can consult texts such as \cite{hunger2007introduction} for an introduction to complex differentials. First, a few definitions:
\begin{itemize}
    \item A matrix Hermitian conjugate  or matrix conjugate transpose is denoted as $\bm{X}^{\dagger}$ and obtained by taking the conjugate of each entry and then transposing the matrix. In the real case the conjugate transpose of matrix $\bm{X}$ is its transpose $\bm{X}^T$.
    \item For a square complex matrix X, we define $symh(\bm{X}) = \frac{\bm{X} + \bm{X}^{\dagger}}{2}$ where $symh(\bm{X})$ is a Hermitian matrix. This is the analogue of the $sym(\cdot)$ operator in the real case.
    \item If a matrix satisfies $\bm{X} = - \bm{X}^{\dagger}$, $\bm{X}$ is called antihermitian.
    \item We define the operator $hcopyltu(\bm{M}) = symh(\bm{M} \circ \bm{E})$, with $\bm{M}$ a complex square matrix and $\bm{E}$ is the mask matrix used in the Hadamard product in the real case and defined in Section 2.3. In practice, $hcopyltu(\cdot)$ conjugates the lower triangle elements of $\bm{M}$ and copies them to the upper triangle. Also, it sets the diagonal of $\bm{M}$ to its real part (sets the imaginary part to zero). This is analogous to the $copyltu(\cdot)$ operator in the real case.
\end{itemize}

Matching the range of our loss, consider a real valued function $f(z)$, where $z \in \mathbb{C}$, $z = x + iy$ and $z^{*} = x - iy$. Under suitable conditions listed in \cite{hunger2007introduction}, we can write:
\begin{equation} \label{wirtinger}
    df = \frac{\partial f}{\partial z}dz + \frac{\partial f}{\partial z^{*}}dz^{*}
\end{equation}
The partial derivatives (Wirtinger derivatives) $\frac{\partial f}{\partial z}$ and $\frac{\partial f}{\partial z^{*}}$ are defined in \cite{hunger2007introduction}. Analogizing that $f$ is our loss $\mathcal{L}$ and $z$ is our complex $\bm{A}$ matrix, the introduction of the $h.c.$ terms in Equation \ref{trid2} produces a real valued trace to match the real valued $d\mathcal{L}$ and the correct gradient $\bar{\bm{A}}$.

\subsection{Deep and Square Complex Input Matrices QR Matrix Backprop}

Our results for complex QR reversed mode auto-diff for square and deep complex matrices agree with the results given in \cite{blog} and \cite{hubig2019use}. We provide the complete matrix backprop process steps derivations in the Appendix.

For the Step 1 of the matrix backpropagation process, similarly to the derivations for the real case in Section 3.1, we get
\begin{equation*}
\begin{split}
    d\bm{Q} &= (d\bm{A} - \bm{Q}d\bm{R})\bm{R}^{-1} \\
    d\bm{R} &= (symh(\bm{C})\circ \bm{E}^{T})\bm{R}
\end{split}
\end{equation*}
with $\bm{C}=\bm{Q}^{\dagger}d\bm{A}\bm{R}^{-1}.$
The additional assumptions required for the complex version of Step 1 are that both $\bm{R}$ and $d\bm{R}$ have real diagonal elements.

For Step 2 of the matrix backpropagation process, the Hermitian conjugates of all the terms in the trace identity must be considered, as we discussed in Section 4.1. Appendix 6.5 presents the full derivations for the reader who wants to delve deeper. The derived gradient expression is
\begin{equation}\label{eqn:deep_dagger}
\begin{split}
\bar{\bm{A}} & = \left[\bar{\bm{Q}} + \bm{Q}symh(\bm{M}\circ \bm{E})\right]\bm{R}^{-\dagger} \\
& = \left[\bar{\bm{Q}} + \bm{Q}hcopyltu(\bm{M})\right]\bm{R}^{-\dagger}.
\end{split}
\end{equation}
where $\bm{M} = \bm{R}\bar{\bm{R}}^{\dagger} - \bar{\bm{Q}}^{\dagger}\bm{Q}$. 

\subsection{Alternative Formulation with Correction for the Square and Deep Complex QR Case}

Next we discuss the correction needed for the Equation 42 given in \cite{walter2010algorithmic} for the real case to be adapted to the complex case. 
In \cite{walter2010algorithmic}, the derivation omits the $h.c.$ terms for the complex case that need to be included in the trace identity to produce the equivalent complex gradient. A full re-derivation would be redundant given derivations already in this paper. However, determining this correction for the complex case is necessary because Equation 42 of \cite{walter2010algorithmic} is implemented by TensorFlow and PyTorch for the real cases. We include this correction as the smallest code change that was necessary to get the correct gradient for the complex QR when using Equation 42 from \cite{walter2010algorithmic}.

To recap, Equation 42 in \cite{walter2010algorithmic} is 
\begin{equation*}
    \bar{\bm{A}} = \bm{Q}(\bar{\bm{R}} + \bm{P}_{L} \circ (\bm{R}\bar{\bm{R}}^{T} - \bar{\bm{R}}\bm{R}^{T} + \bm{Q}^{T} \bar{\bm{Q}} - \bar{\bm{Q}}^{T} \bm{Q})\bm{R}^{-T}) + (\bar{\bm{Q}} - \bm{Q}\bm{Q}^{T} \bar{\bm{Q}})\bm{R}^{-T},
\end{equation*}
with $\bm{P}_L = (i>j)$, a strictly lower tridiagonal matrix with all ones beneath the diagonal and zeroes along and above the main diagonal. Here $\bm{A}$ is square or deep and both real and full rank.

Naively considering the Hermitian transpose in place of the transpose :
\begin{equation*}
    \bar{\bm{A}} = \bm{Q}(\bar{\bm{R}} + \bm{P}_{L} \circ (\bm{R}\bar{\bm{R}}^{\dagger} - \bar{\bm{R}}\bm{R}^{\dagger} + \bm{Q}^{\dagger} \bar{\bm{Q}} - \bar{\bm{Q}}^{\dagger} \bm{Q})\bm{R}^{-\dagger}) + (\bar{\bm{Q}} - \bm{Q}\bm{Q}^{\dagger} \bar{\bm{Q}})\bm{R}^{-\dagger}.
\end{equation*}

We can follow the term rearranging steps outlined for the real case in Appendix 6.2 to get
\begin{equation}\label{parans}
\bar{\bm{A}} = \left[\bar{\bm{Q}} + \bm{Q}(\bm{P}_{L} \circ (\bm{M} - \bm{M}^{\dagger}) + \bm{M}^{\dagger})\right]\bm{R}^{-\dagger}.
\end{equation}
Note that $\bm{M} = \bm{R}\bar{\bm{R}}^{\dagger} - \bar{\bm{Q}}^{\dagger}\bm{Q}$ and $\bm{E} = 2\bm{P}_{L} + \bm{I}$ by definition, so one can write $\bm{P}_{L} = 0.5(\bm{E} - \bm{I})$. Note that the $\bm{Q}$ multiplier in Equation \ref{parans} is not equal to the corresponding term in the correct gradient expression given in Equation \ref{eqn:deep_dagger} in Section 4.3.,  
\begin{equation}\label{rhs}
\begin{split}
\bm{P}_{L} \circ (\bm{M} - \bm{M}^{\dagger}) + \bm{M}^{\dagger} \neq symh(\bm{M}\circ \bm{E}).
\end{split}
\end{equation}
We next derive a correction term $\mathcal{C}$ so that:

\begin{equation}
\mathcal{C} + (\bm{P}_{L} \circ (\bm{M} - \bm{M}^{\dagger}) + \bm{M}^{\dagger}) = symh(\bm{M}\circ \bm{E}).
\end{equation}

So
\begin{equation}\label{eqn:correction}
\mathcal{C} = symh(\bm{M}\circ \bm{E}) - (\bm{P}_{L} \circ (\bm{M} - \bm{M}^{\dagger}) + \bm{M}^{\dagger}).
\end{equation}

First simplifying the $symh(\bm{M}\circ \bm{E})$ term and recalling that $\bm{E} = 2\bm{P}_l + \bm{I}$

\begin{equation*}
\begin{split}
symh(\bm{M}\circ \bm{E}) &=0.5(\bm{M} \circ \bm{E}) + 0.5(\bm{M} \circ \bm{E})^{\dagger} \\
& = \bm{M} \circ \bm{P}_l + 0.5(\bm{M} \circ \bm{I}) + (\bm{M} \circ \bm{P}_l)^{\dagger} + 0.5(\bm{M} \circ \bm{I})^{\dagger}.
\end{split}
\end{equation*}

Continuing the matrix acrobatics, 
\begin{equation*}
\begin{split}
\mathcal{C} &= \bm{M} \circ \bm{P}_l + 0.5(\bm{M} \circ \bm{I}) + (\bm{M} \circ \bm{P}_l)^{\dagger} + 0.5(\bm{M} \circ \bm{I})^{\dagger} - \bm{P}_{L} \circ \bm{M} + \bm{P}_{L} \circ \bm{M}^{\dagger} - \bm{M}^{\dagger} \\
&= 0.5(\bm{M} \circ \bm{I}) + (\bm{M} \circ \bm{P}_l)^{\dagger} + 0.5(\bm{M} \circ \bm{I})^{\dagger} + \bm{P}_{L} \circ \bm{M}^{\dagger} - \bm{M}^{\dagger}.
\end{split}
\end{equation*}

Note that one can write $(\bm{M} \circ \bm{P}_l)^{\dagger} + \bm{P}_{L} \circ \bm{M}^{\dagger} = \bm{M}^{\dagger} - \bm{M}^{\dagger} \circ \bm{I}$. Furthermore, writing $0.5(\bm{M} \circ \bm{I}) + 0.5(\bm{M} \circ \bm{I})^{\dagger} = \Re(\bm{M} \circ \bm{I})$, we simplify the correction to get
\begin{equation*}
\begin{split}
\mathcal{C} &= \bm{M}^{\dagger} - \bm{M}^{\dagger} \circ \bm{I} + \Re(\bm{M} \circ \bm{I}) - \bm{M}^{\dagger} \\
& = \Re(\bm{M} \circ \bm{I}) - \bm{M}^{\dagger} \circ \bm{I} \\
&=  \bm{M} \circ \bm{I} - \Re(\bm{M} \circ \bm{I})\\
&= i \Im(\textnormal{diag}(\bm{M})).
\end{split}
\end{equation*}

In summary, the correction $\mathcal{C}$ is a matrix of purely complex numbers (real parts zero) and the imaginary part of each matrix entry is equal to the imaginary part of the diagonal of $\bm{M}$.

\subsection{Wide Complex Input Matrix QR Backprop}

The matrix backpropagation process for wide complex input matrix QR is similar to the real case derived through the partitioning trick and given in Equation 9 and 10. The complex matrix version is
\begin{equation} \label{eqn:wide_dagger}
\bar{\bm{A}} = \left[(\bar{\bm{Q}}_{prime} + \bm{Q} hcopyltu(\bm{M}))\bm{U}^{-\dagger}|\bm{Q}\bar{\bm{V}}\right],
\end{equation} 
with the partitioning trick $\bm{R}=[\bm{U}|\bm{V}]$, and $\bar{\bm{Q}}_{prime} = \bar{\bm{Q}} + \bm{Y}\bar{\bm{V}}^{\dagger}$, 
and $\bm{M} = \bm{U}\bar{\bm{U}}^{\dagger} - \bar{\bm{Q}}_{prime}^{\dagger}\bm{Q}$. A full rederivation must again consider the extra Hermitian conjugate terms in the trace identity. For brevity, we omit the derivation of Equation \ref{eqn:wide_dagger} but note it follows according to the reasoning used to arrive at Proposition 2. To the best of our knowledge, we are unaware of the wide case complex QR auto-diff derivations elsewhere.

\section{Final Notes}
In this article we presented the reverse mode automatic differentiation algorithms for QR and LQ decompositions (reduced mode) that can be used in matrix \bp\ derivations for any full-rank matrix orders. As a result of this work the QR decomposition (forward and backward, for all input shapes real matrices) was implemented in \href{https://github.com/apache/incubator-mxnet/blob/master/src/operator/numpy/linalg/np_qr-inl.h}{MXNet} \cite{chen2015mxnet}. 
In TensorFlow Core \cite{abadi2016tensorflow} and PyTorch \cite{paszke2019pytorch} the gradient for wide inputs and complex matrices of all shapes were implemented.
All frameworks support CPU and GPU implementations for the real case. 
In these frameworks the QR decomposition can be applied in batches and performed on the last two dimensions of larger tensors. The software implementation of the methods across the popular deep learning frameworks allows researchers and engineers to use the differentiable QR across common frameworks.

\bibliographystyle{plain}
\bibliography{qr.bib}

\begin{thebibliography}{10}

\bibitem{abadi2016tensorflow}
Mart{\'\i}n Abadi, Paul Barham, Jianmin Chen, Zhifeng Chen, Andy Davis, Jeffrey
  Dean, Matthieu Devin, Sanjay Ghemawat, Geoffrey Irving, Michael Isard, et~al.
\newblock Tensorflow: A system for large-scale machine learning.
\newblock In {\em 12th $\{$USENIX$\}$ Symposium on Operating Systems Design and
  Implementation ($\{$OSDI$\}$ 16)}, pages 265--283, 2016.

\bibitem{anderson1999lapack}
Edward Anderson, Zhaojun Bai, Christian Bischof, Susan Blackford, Jack
  Dongarra, Jeremy Du~Croz, Anne Greenbaum, Sven Hammarling, Alan McKenney, and
  Danny Sorensen.
\newblock {\em {LAPACK} Users' guide}, volume~9.
\newblock Siam, 1999.

\bibitem{chen2015mxnet}
Tianqi Chen, Mu~Li, Yutian Li, Min Lin, Naiyan Wang, Minjie Wang, Tianjun Xiao,
  Bing Xu, Chiyuan Zhang, and Zheng Zhang.
\newblock {MXN}et: A flexible and efficient machine learning library for
  heterogeneous distributed systems.
\newblock {\em arXiv preprint arXiv:1512.01274}, 2015.

\bibitem{9026903}
Z.~{Dang}, K.~M. {Yi}, Y.~{Hu}, F.~{Wang}, P.~{Fua}, and M.~{Salzmann}.
\newblock Eigendecomposition-free training of deep networks for linear
  least-square problems.
\newblock {\em IEEE Transactions on Pattern Analysis and Machine Intelligence},
  pages 1--1, 2020.

\bibitem{el2020orthonet}
Mireille El~Gheche, Giovanni Chierchia, and Pascal Frossard.
\newblock Orthonet: Multilayer network data clustering.
\newblock {\em IEEE Transactions on Signal and Information Processing over
  Networks}, 6:13--23, 2020.

\bibitem{giles2008extended}
Mike~B Giles.
\newblock Collected matrix derivative results for forward and reverse mode
  algorithmic differentiation.
\newblock In {\em Advances in Automatic Differentiation}, pages 35--44.
  Springer, 2008.

\bibitem{hubig2019use}
Claudius Hubig.
\newblock Use and implementation of autodifferentiation in tensor network
  methods with complex scalars.
\newblock {\em arXiv preprint arXiv:1907.13422}, 2019.

\bibitem{hunger2007introduction}
Raphael Hunger.
\newblock An introduction to complex differentials and complex
  differentiability.
\newblock 2007.

\bibitem{ionescu2015training}
Catalin Ionescu, Orestis Vantzos, and Cristian Sminchisescu.
\newblock Training deep networks with structured layers by matrix
  backpropagation.
\newblock {\em arXiv preprint arXiv:1509.07838}, 2015.

\bibitem{kanakis2020reparameterizing}
Menelaos Kanakis, David Bruggemann, Suman Saha, Stamatios Georgoulis, Anton
  Obukhov, and Luc Van~Gool.
\newblock Reparameterizing convolutions for incremental multi-task learning
  without task interference.
\newblock In {\em European Conference on Computer Vision}, pages 689--707.
  Springer, 2020.

\bibitem{li2017reconstruction}
Lixiang Li, Dafei Xu, Haipeng Peng, J{\"u}rgen Kurths, and Yixian Yang.
\newblock Reconstruction of complex network based on the noise via qr
  decomposition and compressed sensing.
\newblock {\em Scientific reports}, 7(1):1--13, 2017.

\bibitem{liao2019differentiable}
Hai-Jun Liao, Jin-Guo Liu, Lei Wang, and Tao Xiang.
\newblock Differentiable programming tensor networks.
\newblock {\em Physical Review X}, 9(3):031041, 2019.

\bibitem{blog}
Jin-Guo Liu.
\newblock Linear algebra autodiff (complex valued).
\newblock \url{https://giggleliu.github.io/2019/04/02/einsumbp.html}.
\newblock Accessed: 2020-12-09.

\bibitem{ma2019graph}
Yao Ma, Suhang Wang, Charu~C Aggarwal, and Jiliang Tang.
\newblock Graph convolutional networks with eigenpooling.
\newblock In {\em Proceedings of the 25th ACM SIGKDD International Conference
  on Knowledge Discovery \& Data Mining}, pages 723--731, 2019.

\bibitem{murphy2012machine}
Kevin~P Murphy.
\newblock {\em Machine learning: a probabilistic perspective}.
\newblock MIT press, 2012.

\bibitem{nvidia2011nvidia}
CUDA Nvidia.
\newblock Nvidia cuda c programming guide.
\newblock {\em Nvidia Corporation}, 120(18):8, 2011.

\bibitem{oliphant2006guide}
Travis~E Oliphant.
\newblock {\em A guide to {N}um{P}y}, volume~1.
\newblock Trelgol Publishing USA, 2006.

\bibitem{olteanu2008technical}
Denisa Olteanu and Laura~J Freeman.
\newblock Technical report on the evaluation of median rank regression and
  maximum likelihood estimation techniques for a two-parameter {Weibull}
  distribution.
\newblock Technical report, Virginia Tech, 2008.

\bibitem{pan2009survey}
Sinno~Jialin Pan and Qiang Yang.
\newblock A survey on transfer learning.
\newblock {\em IEEE Transactions on knowledge and data engineering},
  22(10):1345--1359, 2009.

\bibitem{paszke2019pytorch}
Adam Paszke, Sam Gross, Francisco Massa, Adam Lerer, James Bradbury, Gregory
  Chanan, Trevor Killeen, Zeming Lin, Natalia Gimelshein, Luca Antiga, et~al.
\newblock Pytorch: An imperative style, high-performance deep learning library.
\newblock In {\em Advances in neural information processing systems}, pages
  8026--8037, 2019.

\bibitem{ramon2004numeric}
Jan Ramon and Kurt Driessens.
\newblock On the numeric stability of gaussian processes regression for
  relational reinforcement learning.
\newblock In {\em ICML-2004 Workshop on Relational Reinforcement Learning},
  pages 10--14, 2004.

\bibitem{searle2017matrix}
Shayle~R Searle and Andre~I Khuri.
\newblock {\em Matrix algebra useful for statistics}.
\newblock John Wiley \& Sons, 2017.

\bibitem{seeger2017auto}
Matthias Seeger, Asmus Hetzel, Zhenwen Dai, Eric Meissner, and Neil~D Lawrence.
\newblock Auto-differentiating linear algebra.
\newblock {\em arXiv preprint arXiv:1710.08717}, 2017.

\bibitem{trefethen1997numerical}
Lloyd~N Trefethen and David Bau~III.
\newblock {\em Numerical linear algebra}, volume~50.
\newblock Siam, 1997.

\bibitem{van1983matrix}
Charles~F Van~Loan and Gene~H Golub.
\newblock {\em Matrix computations}.
\newblock Johns Hopkins University Press Baltimore, 1983.

\bibitem{virtanen2020scipy}
Pauli Virtanen, Ralf Gommers, Travis~E Oliphant, Matt Haberland, Tyler Reddy,
  David Cournapeau, Evgeni Burovski, Pearu Peterson, Warren Weckesser, Jonathan
  Bright, et~al.
\newblock Scipy 1.0: fundamental algorithms for scientific computing in python.
\newblock {\em Nature methods}, 17(3):261--272, 2020.

\bibitem{walter2010algorithmic}
Sebastian~F Walter and Lutz Lehmann.
\newblock Algorithmic differentiation of linear algebra functions with
  application in optimum experimental design (extended version).
\newblock {\em arXiv preprint arXiv:1001.1654}, 2010.

\bibitem{yasotharan2010simple}
Ambighairajah Yasotharan.
\newblock A simple expression for the matrix gradient of a diagonal element of
  r in qr decomposition for use in mimo communications and signal processing.
\newblock Technical report, DEFENCE RESEARCH AND DEVELOPMENT CANADA OTTAWA
  (ONTARIO), 2010.

\end{thebibliography}

\clearpage

\section {Appendix}

\subsection{Lemma 1.}
For a deep matrix with QR decomposition $\bm{A}=\bm{Q}\bm{R}$, the variations of $\bm{Q}$ and $\bm{R}$ for a given variation $d\bm{A}$ of $\bm{A}$, denoted $d\bm{Q}$ and $d\bm{R}$ respectively,  are 
\begin{equation}
\begin{split}
d\bm{Q} &= (d\bm{A} - \bm{Q}d\bm{R})\bm{R}^{-1} \\
d\bm{R} &= (sym(\bm{C})\circ \bm{E}^T)\bm{R},
\end{split}
\end{equation}
with $\bm{C}=\bm{Q}^{T}d\bm{A}\bm{R}^{-1}.$

\textit{Proof}
First calculate the variations of $\bm{Q}$ and $\bm{R}$, where the following hold:
\begin{itemize}
\item $\bm{R}$ and $d\bm{R}$ are square, upper triangular and full rank.
\item $\bm{Q}$ is orthogonal with $\bm{Q}^T\bm{Q} = \bm{I}$ and hence (after taking the first variation)
\begin{equation}\label{eq:13}
    d\bm{Q}^T\bm{Q} +\bm{Q}^Td\bm{Q} = \bm{0}
\end{equation}
so $\bm{Q}^{T}d\bm{Q}$ is skew-symmetric with $\bm{Q}^{T}d\bm{Q} = -d\bm{Q}^T\bm{Q}$.
\end{itemize}

Take the first variation of the QR decomposition $\bm{A} = \bm{Q}\bm{R}$
\begin{equation*}
    d\bm{A} = d\bm{Q}\bm{R} + \bm{Q}d\bm{R}.
\end{equation*}

Right multiply by $\bm{R}^{-1}$ to get
\begin{equation}\label{eq:skew}
    d\bm{Q} = (d\bm{A} - \bm{Q}d\bm{R})\bm{R}^{-1}.
\end{equation}

Next, find an expression for $d\bm{R}$. First left-multiply Equation \ref{eq:skew} by $\bm{Q}^{T}$ to get
\begin{equation*}
    \bm{Q}^{T}d\bm{Q} = \bm{Q}^{T}(d\bm{A} - \bm{Q}d\bm{R})\bm{R}^{-1}.
\end{equation*}
Now $\bm{Q}^{T}d\bm{Q}$ is skew-symmetric, so the \rhs\ is also skew-symmetric, and
\begin{equation*}
    \bm{Q}^{T}(d\bm{A} - \bm{Q}d\bm{R})\bm{R}^{-1} = - \bm{R}^{-T}(d\bm{A} - \bm{Q}d\bm{R})^{T}\bm{Q}
\end{equation*}
which simplifies to 
\begin{equation*}
    \bm{Q}^{T}d\bm{A}\bm{R}^{-1} + (\bm{Q}^{T}d\bm{A}\bm{R}^{-1})^{T} = d\bm{R}\bm{R}^{-1} + {(d\bm{R}\bm{R}^{-1})}^T.
\end{equation*}
Consequently, with some left to right algebraic manipulations,
\begin{equation}\label{eq:14}
     sym(\bm{C}) = sym(d\bm{R}\bm{R}^{-1}),
\end{equation}
with $\bm{C}=\bm{Q}^{T}d\bm{A}\bm{R}^{-1}.$ In Equation \ref{eq:14}, $d\bm{R}\bm{R}^{-1}$ is upper triangular so we express $d\bm{R}$, using the notation of Lemma 1, as $d\bm{R} = (sym(\bm{C})\circ \bm{E}^{T})\bm{R}.$

\subsection{Equivalence of Equation \ref{eqn:a_bar} in Proposition 1 and Equation \ref{eqn:walter} in Section 3.1}

In Section 3.1 we derived the gradient for the QR decomposition of square and deep input matrices and noted that the resulting formula is equivalent to the one given in \ref{eqn:walter}, albeit more compact and computationally efficient. We prove the equivalence next.

\textit{Proof}

We need to prove that Equation \ref{eqn:walter},
\begin{equation*}
    \bar{\bm{A}} = \bm{Q}(\bar{\bm{R}} + \bm{P}_{L} \circ (\bm{R}\bar{\bm{R}}^{T} - \bar{\bm{R}}\bm{R}^{T} + \bm{Q}^{T} \bar{\bm{Q}} - \bar{\bm{Q}}^{T} \bm{Q})\bm{R}^{-T}) + (\bar{\bm{Q}} - \bm{Q}\bm{Q}^{T} \bar{\bm{Q}})\bm{R}^{-T},
\end{equation*}
with $\bm{P}_L = (i>j)$, a strictly lower tridiagonal matrix with all ones below the diagonal and zeroes along and above the main diagonal, simplifies to Equation \ref{eqn:a_bar},
\begin{equation*}
\bar{\bm{A}} = \left[\bar{\bm{Q}} + \bm{Q}copyltu(\bm{M})\right]\bm{R}^{-T},   
\end{equation*}
where $\bm{M} = \bm{R}\bar{\bm{R}}^T - \bar{\bm{Q}}^{T}\bm{Q}$ for real matrices (it does not for complex matrices).

First, in Equation \ref{eqn:walter}, we rearrange and simplify the innermost parenthesis as
\begin{equation*}
\begin{split}
    &=\bm{R}\bar{\bm{R}}^{T} - \bar{\bm{R}}\bm{R}^{T} + \bm{Q}^{T} \bar{\bm{Q}} - \bar{\bm{Q}}^{T} \bm{Q} \\
    &= (\bm{R}\bar{\bm{R}}^{T} - \bar{\bm{Q}}^{T} \bm{Q}) - (\bar{\bm{R}}\bm{R}^{T} - \bm{Q}^{T} \bar{\bm{Q}}) \\
    &= \bm{M} - \bm{M}^{T}.
\end{split}
\end{equation*}

Now we can re-write Equation \ref{eqn:walter} as 
\begin{equation*}
 \bar{\bm{A}} = \bm{Q}(\bar{\bm{R}}\bm{R}^T + \bm{P}_{L} \circ (\bm{M} - \bm{M}^{T})\bm{R}^{-T}) + (\bar{\bm{Q}} - \bm{Q}\bm{Q}^{T} \bar{\bm{Q}})\bm{R}^{-T}.
 \end{equation*}

Next we factor out $\bm{R}^{-T}$ 
\begin{equation*}
 \bar{\bm{A}} = [\bm{Q}(\bar{\bm{R}}\bm{R}^{T} + \bm{P}_{L} \circ ( \bm{M} - \bm{M}^{T})) + \bar{\bm{Q}} - \bm{Q}\bm{Q}^{T} \bar{\bm{Q}}]\bm{R}^{-T}.
 \end{equation*}
 Then we re-arrange terms and factor out $\bm{Q}$
 \begin{equation*}
 \bar{\bm{A}} = [\bar{\bm{Q}} + \bm{Q}( \bm{P}_{L} \circ ( \bm{M} - \bm{M}^{T}) + \bar{\bm{R}}\bm{R}^{T} - \bm{Q}^{T} \bar{\bm{Q}})]\bm{R}^{-T}.
 \end{equation*}

We now identify matrix $\bm{M}^T$ and rewrite 
\begin{equation*}
\bar{\bm{A}} = \left[\bar{\bm{Q}} + \bm{Q}(\bm{P}_{L} \circ (\bm{M} - \bm{M}^{T}) + \bm{M}^{T})\right]\bm{R}^{-T}.
\end{equation*}
Note that $\bm{E} = 2\bm{P}_{L} + \bm{I}$ by definition, so one can write $\bm{P}_{L} = 0.5(\bm{E} - \bm{I})$ and 

\begin{equation*}
\begin{split}
\bar{A} &= \left[\bar{\bm{Q}} + \bm{Q}\left(0.5(\bm{E} - \bm{I}) \circ \bm{M} - 0.5(\bm{E} -\bm{I}) \circ \bm{M}^T + \bm{M}^{T}\right)\right]\bm{R}^{-T} \\
&=\left[\bar{\bm{Q}} + \bm{Q}\left(0.5(\bm{M} \circ E) - 0.5(\bm{I} \circ \bm{M}) - 0.5(\bm{M}^T \circ \bm{E}) + 0.5(\bm{I} \circ \bm{M}^T) + \bm{M}^T\right)\right]\bm{R}^{-T} \\
&= \left[\bar{\bm{Q}} + \bm{Q}\left(0.5(\bm{E} \circ \bm{M} + (\bm{E} \circ \bm{M})^{T}\right)\right]\bm{R}^{-T} \\
&= \left[\bar{\bm{Q}} + \bm{Q}sym(\bm{E}\circ \bm{M})\right]\bm{R}^{-T} \\
&= \left[\bar{\bm{Q}} + \bm{Q}copyltu(\bm{M})\right]\bm{R}^{-T},
\end{split}
\end{equation*}
which proves the equivalence of Equation \ref{eqn:a_bar} and \ref{eqn:walter} in the real input matrices case. The equivalence does not hold in the complex case we discuss a correction in Section 4.

\subsection{Lemma 2.}
When $\bm{A}$ is wide, full rank, and with the first $k$ columns forming a square full rank matrix, the variations of $\bm{A}$ and $\bm{R}$ can be partitioned similarly to $\bm{A}$ and $\bm{R}$ as $d\bm{A}=[d\bm{X}|d\bm{Y}]$ and $d\bm{R} = [d\bm{U}|d\bm{V}]$ with corresponding orders such that 
\begin{equation*}
\begin{split}
    d\bm{Q} = (d\bm{X} - \bm{Q}d\bm{U})\bm{U}^{-1}, \\ 
    d\bm{V} = \bm{Q}^{T}(d\bm{Y} - d\bm{Q}\bm{V}), \\
    d\bm{U} = (sym(\bm{C})\circ \bm{E}^{T})\bm{U}, \\
\end{split}
\end{equation*}
with $\bm{C}=\bm{Q}^{T}d\bm{X}\bm{U}^{-1}.$

\textit{Proof.}
We want to calculate the variations $d\bm{Q}$ and $[d\bm{U}|d\bm{V}]$ for given $[d\bm{X}|d\bm{Y}]$, where

\begin{itemize}
    \item $d\bm{R}$ and $d\bm{U}$ are upper triangular matrices.
    \item $\bm{Q}$ is orthogonal with
    \begin{equation}\label{eq:11}
        d\bm{Q}^{T}\bm{Q} + \bm{Q}^{T}d\bm{Q} = \bm{0}.
    \end{equation}
\end{itemize}
The exposition is continued using partitioned matrices. 
Taking the first variation of the QR decomposition of $\bm{X}$, we have 
\begin{equation*}
    d\bm{X} = d\bm{Q}\bm{U} + \bm{Q}d\bm{U}.
\end{equation*}
Then right multiply by $\bm{U}^{-1}$ to get
\begin{equation}\label{dU}
    d\bm{Q} = (d\bm{X} - \bm{Q}d\bm{U})\bm{U}^{-1}.
\end{equation}

Similarly take the first variation of $\bm{Y}=\bm{Q}\bm{V}$ in Equation \ref{eqn:prop2} to get $d\bm{Y} = \bm{Q}d\bm{V} + d\bm{Q}\bm{V}$ and then $d\bm{V} = \bm{Q}^{T}(d\bm{Y} - d\bm{Q}\bm{V}).$
Next, we use Equation \ref{dU} to find an expression for $d\bm{U}$. First left-multiply by $\bm{Q}^T$ to get
\begin{equation*}
    \bm{Q}^{T}d\bm{Q} = \bm{Q}^{T}(d\bm{X} - \bm{Q}d\bm{U})\bm{U}^{-1}.
\end{equation*}
Now the \lhs\ is skew-symmetric, so the \rhs\ is also skew-symmetric, with
\begin{equation*}
    \bm{Q}^{T}(d\bm{X} - \bm{Q}d\bm{U})\bm{U}^{-1} = - \bm{U}^{-T} (d\bm{X} - \bm{Q}d\bm{U})^{T}\bm{Q},
\end{equation*}
which simplifies to 
\begin{equation*}
    \bm{Q}^{T}d\bm{X}\bm{U}^{-1} + (\bm{Q}^{T}d\bm{X}\bm{U}^{-1})^{T} = d\bm{U}\bm{U}^{-1} + {(d\bm{U}\bm{U}^{-1})}^T.
\end{equation*}
Then write
\begin{equation*}
    sym(d\bm{U}\bm{U}^{-1}) = sym(\bm{C}),
\end{equation*}
with $\bm{C}=\bm{Q}^{T}d\bm{X}\bm{U}^{-1}.$ Matrix $d\bm{U}\bm{U}^{-1}$ is upper triangular so we express $d\bm{U}$, using the notation in Lemma 2,
\begin{equation*}
    d\bm{U} = (sym(\bm{C})\circ \bm{E}^T)\bm{U}.
\end{equation*}

\subsection{Proposition 3. LQ BP For Deep Input Matrices}

\textit{Proof.}
The proof follows closely the proof of Proposition 2.

If $\bm{A}$ is deep, of rank $k$, with the top $k$ rows forming a square full-rank matrix, and partitioned as illustrated in Section 3.3, Proposition 3, then
\begin{equation}\label{eq:9}
\bm{X} = \bm{U}\bm{Q}
\end{equation}
is the LQ decomposition of $\bm{X}$, with $\bm{Q}$ and $\bm{U}$ square, full-rank invertible matrices. Then $\bm{Y} = \bm{V}\bm{Q}$. Hence we can apply the chain rule and the two step matrix \bp\ to obtain $\bar{\bm{A}}$.

\textbf{BP Step 1: Variations} 
First calculate the variations of $\bm{L}$ and $\bm{Q}$, $d\bm{Q}$ and $d\bm{L}$ for a given variation $d\bm{A}$ of $\bm{A}$. The variations of $\bm{A}$ and $\bm{L}$ can also be partitioned as $d\bm{A}=[d\bm{X}|d\bm{Y}]$ and $d\bm{L} = [d\bm{U}|d\bm{V}]$ with corresponding orders such that
\begin{equation*}
    d\bm{Q} = \bm{U}^{-1}(d\bm{X} - d\bm{U}\bm{Q})
\end{equation*}
\begin{equation*}
    d\bm{V} = (d\bm{Y} - \bm{V}d\bm{Q})\bm{Q}^{T}.
\end{equation*}
\begin{equation*}
    d\bm{U} = \bm{U}(sym(\bm{C})\circ \bm{E}).
\end{equation*}
with $\bm{C}=\bm{U}^{-1}d\bm{X}\bm{Q}^{T}.$

\textit{Proof.}
We want to calculate the variations $d\bm{Q}$ and $[d\bm{U}$, and $d\bm{V}]$ for given $[d\bm{X}$ and $d\bm{Y}]$, where

\begin{itemize}
    \item $d\bm{L}$ and $d\bm{U}$ are lower triangular matrices.
    \item $\bm{Q}$ is orthogonal with
    \begin{equation}\label{eqn:11}
        d\bm{Q}\bm{Q}^T + \bm{Q}d\bm{Q}^T = \bm{0}
    \end{equation}
\end{itemize}

Taking the first variation of the LQ decomposition of $\bm{X}$, we have 
\begin{equation}\label{eqn:one}
    d\bm{X} = d\bm{U}\bm{Q} + \bm{U}d\bm{Q}.
\end{equation}
Then left multiply by $\bm{U}^{-1}$ to get
\begin{equation*}
    d\bm{Q} = \bm{U}^{-1}(d\bm{X} - d\bm{U}\bm{Q})
\end{equation*}

Similarly take the first variation of $\bm{Y}=\bm{V}\bm{Q}$ to get $d\bm{Y} = d\bm{V}\bm{Q} + \bm{V}d\bm{Q}$ and then $d\bm{V} = (d\bm{Y} - \bm{V}d\bm{Q})\bm{Q}^{T}.$
Next, find an expression for $d\bm{U}$. First right-multiply by $\bm{Q}^T$ in Equation \ref{eqn:one} to get
\begin{equation*}
    d\bm{Q}\bm{Q}^{T} = \bm{U}^{-1}(d\bm{X} - \bm{Q}d\bm{U})\bm{Q}^{T}.
\end{equation*}
The \lhs\ is skew-symmetric, so the \rhs\ is also skew-symmetric, and we can write
\begin{equation*}
    \bm{U}^{-1} d\bm{X}\bm{Q}^{T} + (\bm{U}^{-1} d\bm{X}\bm{Q}^{T})^{T} = \bm{U}^{-1}d\bm{U} + {(\bm{U}^{-1}d\bm{U})}^T.
\end{equation*}
Then
\begin{equation*}
    sym(\bm{U}^{-1}d\bm{U}) = sym(\bm{C}),
\end{equation*}
with $\bm{C}=\bm{U}^{-1}d\bm{X}\bm{Q}^{T}.$ Matrix $\bm{U}^{-1}d\bm{U}$ is lower triangular so we can express $d\bm{U}$
\begin{equation*}
    d\bm{U} = \bm{U}(sym(\bm{C})\circ \bm{E}).
\end{equation*}

\textbf{BP Step 2: Partial Derivatives}
We proceed to the second step of the process for matrix \bp\ process to get $\bar{\bm{A}}$. 
We use the trace identity $Tr(\bar{\bm{A}}^{T}d\bm{A}) = Tr(\bar{\bm{L}}^{T}d\bm{L}) + Tr(\bar{\bm{Q}}^{T}d\bm{Q})$. 
Replace the partitioned matrices and the calculated variations $d\bm{Q}$, $d\bm{U}$ and $d\bm{V}$ and get the \lhs\ and \rhs\ coefficients of $d\bm{X}$ and $d\bm{Y}$ to find $\bar{\bm{X}}$ and $\bar{\bm{Y}}$.

First replace $d\bm{V}$ to get
\begin{equation*}
Tr(\bar{\bm{Q}}^Td\bm{Q}) + Tr(\bar{\bm{U}}^Td\bm{U}) - Tr(\bar{\bm{V}}^{T}\bm{V}d\bm{Q}\bm{Q}^{T}) + Tr(\bar{\bm{V}}^{T}d\bm{Y}\bm{Q}^{T}).
\end{equation*}

Identify the $d\bm{Y}$ coefficient and get $\bar{\bm{Y}} = \bar{\bm{V}}\bm{Q}$. Next identify $d\bm{\bar{X}}$. We have
\begin{equation*}
Tr(\bar{\bm{X}}^Td\bm{X}) = Tr(\bar{\bm{Q}}^Td\bm{Q}) + Tr(\bar{\bm{U}}^Td\bm{U}) - Tr(\bar{\bm{V}}^{T}\bm{V}d\bm{Q}\bm{Q}^{T}).
\end{equation*}

The last trace can be simplified by replacing $d\bm{Q}\bm{Q}^T$, using skew-symmetry and $\bm{V} = \bm{Y}\bm{Q}^{T}$ from the partitioning argument, the orthogonality of $\bm{Q}$, and ICP and IT properties to get
\begin{align*}
Tr(\bar{\bm{V}}^{T}\bm{V}d\bm{Q}\bm{Q}^{T}) &= -Tr(\bar{\bm{V}}^{T}\bm{Y}\bm{Q}^{T}\bm{Q}d(\bm{Q}^T)) \\
&= - Tr(\bm{Y}^T\bar{\bm{V}}d\bm{Q}). 
\end{align*}

With this result simplify the \rhs\ of the trace equality
\begin{align*} 
Tr(\bar{\bm{Q}}^Td\bm{Q}) + Tr(\bar{\bm{U}}^Td\bm{U}) + Tr(\bm{Y}^T\bar{\bm{V}}d\bm{Q}) \\
& = Tr(\bar{\bm{Q}}^T_{prime}d\bm{Q}) + Tr(\bar{\bm{U}}^Td\bm{U}).
\end{align*}

Since $\bm{X}$, $\bm{U}$ are square, the rest of the argument for the $d\bm{X}$ term is the argument for LQ decomposition gradient with square input and we get $\bar{\bm{X}}
= \bm{U}^{-T}[\bar{\bm{Q}}_{prime} + copyltu(\bm{M})\bm{Q}] $, 
with 
$\bm{M} = \bm{U}^T\bar{\bm{U}} - \bar{\bm{Q}}_{prime}\bm{Q}^T$, and $\bar{\bm{Q}}_{prime} = \bar{\bm{Q}} + \bm{\bar{V}^T}\bm{Y}$ as proposed.

\subsection{BP Process For Square and Deep Complex Input Matrices QR}
\textbf{BP Step 1: Variations} For a deep complex matrix with QR decomposition $\bm{A}=\bm{Q}\bm{R}$, the variations of $\bm{Q}$ and $\bm{R}$ for a given variation $d\bm{A}$ of $\bm{A}$, denoted $d\bm{Q}$ and $d\bm{R}$ respectively,  are 
\begin{equation}
\begin{split}
d\bm{Q} &= (d\bm{A} - \bm{Q}d\bm{R})\bm{R}^{-1} \\
d\bm{R} &= (symh(C)\circ \bm{E}^{\dagger})\bm{R},
\end{split}
\end{equation}
with $\bm{C}=\bm{Q}^{\dagger}d\bm{A}\bm{R}^{-1}.$

\textit{Proof}
First calculate the variations of $\bm{Q}$ and $\bm{R}$, where the following hold:
\begin{itemize}
\item $\bm{R}$ and $d\bm{R}$ are square, upper triangular and full rank with real diagonal.
\item $\bm{Q}$ is orthogonal with $\bm{Q}^{\dagger}\bm{Q} = \bm{I}$ and hence (after taking the first variation)
\begin{equation}
    d\bm{Q}^{\dagger}\bm{Q} +\bm{Q}^{\dagger}d\bm{Q} = \bm{0}
\end{equation}
so $\bm{Q}^{\dagger}d\bm{Q}$ is antihermitian with $\bm{Q}^{\dagger}d\bm{Q} = -d\bm{Q}^{\dagger}\bm{Q}$.
\end{itemize}

Take the first variation of the QR decomposition $\bm{A} = \bm{Q}\bm{R}$,
\begin{equation*}
    d\bm{A} = d\bm{Q}\bm{R} + \bm{Q}d\bm{R}.
\end{equation*}

Right multiply by $\bm{R}^{-1}$ to get
\begin{equation}\label{eq:skew1}
    d\bm{Q} = (d\bm{A} - \bm{Q}d\bm{R})\bm{R}^{-1}.
\end{equation}

Next, find an expression for $d\bm{R}$. First left-multiply Equation \ref{eq:skew1} by $\bm{Q}^{\dagger}$ to get
\begin{equation*}
    \bm{Q}^{\dagger}d\bm{Q} = \bm{Q}^{\dagger}(d\bm{A} - \bm{Q}d\bm{R})\bm{R}^{-1}.
\end{equation*}
Now $\bm{Q}^{\dagger}d\bm{Q}$ is antihermitian, so the \rhs\ is also antihermitian, and
\begin{equation*}
    \bm{Q}^{\dagger}(d\bm{A} - \bm{Q}d\bm{R})\bm{R}^{-1} = - \bm{R}^{-\dagger}(d\bm{A} - \bm{Q}d\bm{R})^{\dagger}\bm{Q}
\end{equation*}
which simplifies to 
\begin{equation*}
    \bm{Q}^{\dagger}d\bm{A}\bm{R}^{-1} + (\bm{Q}^{\dagger}d\bm{A}\bm{R}^{-1})^{\dagger} = d\bm{R}\bm{R}^{-1} + {(d\bm{R}\bm{R}^{-1})}^{\dagger}.
\end{equation*}
Consequently, with some left to right algebraic manipulations,
\begin{equation}\label{eq:141}
     symh(\bm{C}) = symh(d\bm{R}\bm{R}^{-1}),
\end{equation}
with $\bm{C}=\bm{Q}^{\dagger}d\bm{A}\bm{R}^{-1}.$ Matrix $d\bm{R}\bm{R}^{-1}$ is upper triangular with real diagonal so we express $d\bm{R}$, as $d\bm{R} = (symh(\bm{C}))\circ \bm{E}^{T})\bm{R}.$

\textbf{BP Step 2: Partial Derivatives}
In the second step of the process for matrix \bp\ laid out in Section ~\ref{ad} one obtains partial derivatives using the variations from step one and the trace identity \ref{trid1}. The goal is to express $d\bm{Q}$ and $d\bm{R}$ variations on the \rhs\ of the trace identity as a function of $d\bm{A}$ and then identify the coefficients of $d\bm{A}$. We do the same thing in the complex case, keeping track of the extra h.c. terms. First, replace the $d\bm{Q}$ from step one on the \rhs\ of the trace identity in Equation \ref{trid2} to get
\begin{equation*}
 Tr(\bar{\bm{Q}}^{\dagger}d\bm{Q} + \bar{\bm{R}}^{\dagger}d\bm{R} + h.c.) = Tr(\bar{\bm{Q}}^{\dagger}(d\bm{A} - \bm{Q}d\bm{R})\bm{R}^{-1} + \bar{\bm{R}}^{\dagger}d\bm{R} + h.c.).
\end{equation*}
Next replace $d\bm{R}$ from step one, use the ICP and Hadamard properties of trace, and arrange the terms to isolate $d\bm{A}$. Considering the non $d\bm{A}$ terms only,
\begin{equation*}
\begin{split}
 Tr((\bar{\bm{R}}^{\dagger} - \bm{R}^{-1}\bar{\bm{Q}}^{\dagger}\bm{Q})d\bm{R} + h.c)
    &= Tr(\bm{R}^{-1}(\bm{R}\bar{\bm{R}}^{\dagger} - \bar{\bm{Q}}^{\dagger}\bm{Q})d\bm{R} + h.c.) \\
    & = Tr(\bm{R}^{-1}\bm{M}(symh(\bm{C})\circ \bm{E}^{\bm{T}})\bm{R} + h.c.) \\
    &= Tr(\bm{M}(0.5(\bm{C} + \bm{C}^{\dagger}))\circ \bm{E}^{T} + h.c.) \\
    &= Tr(((\bm{M} \circ \bm{E})(0.5(\bm{C} + \bm{C}^{\dagger})) + (\bm{M}\circ \bm{E})^{\dagger}(0.5(\bm{C}+\bm{C}^{\dagger}))) \\
    &= Tr(0.5((\bm{M} \circ E) + (\bm{M}\circ \bm{E})^{\dagger})\bm{C} + h.c.) \\
    &= Tr(0.5((\bm{M} \circ \bm{E}) + h.c.)\bm{Q}^{\dagger}d\bm{A}\bm{R}^{-1} + h.c.)
\end{split}
\end{equation*}

where we define $\bm{M} = \bm{R}\bar{\bm{R}}^{\dagger} - \bar{\bm{Q}}^{\dagger}\bm{Q}$. 
Returning to the full trace identity Equation \ref{trid2}, and focusing only on  terms that multiply $d\bm{A}$ (no h.c. remaining terms) we identify
\begin{equation*}
 Tr(\bar{\bm{A}}^{\dagger}d\bm{A}) = Tr(\bm{R}^{-1}\bar{\bm{Q}}^{\dagger}d\bm{A} + \bm{R}^{-1}(0.5((\bm{M}\circ \bm{E}) + h.c.))\bm{Q}^{\dagger}d\bm{A}).
\end{equation*}
The matrix that left multiplies $d\bm{A}$ is then
\begin{equation*}
\begin{split}
\bar{\bm{A}} & = \left[\bar{\bm{Q}} + \bm{Q}symh(\bm{M}\circ \bm{E})\right]\bm{R}^{-\dagger} \\
& = \left[\bar{\bm{Q}} + \bm{Q}hcopyltu(\bm{M})\right]\bm{R}^{-\dagger}.
\end{split}
\end{equation*}

\end{document}